\documentclass[a4paper]{article}

\usepackage[utf8]{inputenc}
\usepackage[textwidth=16cm]{geometry}
\usepackage{lmodern}
\usepackage{array}
\usepackage{graphicx}
\usepackage{amssymb, amsfonts, amsthm, amsmath}
\usepackage[english]{babel}
\usepackage[pdfborder={0 0 0}]{hyperref}
\usepackage{tikz}
\usetikzlibrary{decorations}
\usetikzlibrary{decorations.pathreplacing}
\tikzset{individu/.style={draw,thick}}

\usepackage{subfigure}

\theoremstyle{plain}
\newtheorem{theorem}{Theorem}[section]
\newtheorem{corollary}[theorem]{Corollary}
\newtheorem{lemma}[theorem]{Lemma}
\newtheorem{proposition}[theorem]{Proposition}

\theoremstyle{definition}

\theoremstyle{remark}
\newtheorem{remark}[theorem]{Remark}

\newcommand{\N}{\mathbb{N}}
\newcommand{\Z}{\mathbb{Z}}
\newcommand{\R}{\mathbb{R}}

\newcommand{\calR}{\mathcal{R}}

\newcommand{\ind}[1]{\mathbf{1}_{\{#1\}}}
\newcommand{\indset}[1]{\mathbf{1}_{#1}}
\newcommand{\floor}[1]{{\left\lfloor #1 \right\rfloor}}

\newcommand{\calC}{\mathcal{C}}

\numberwithin{equation}{section}

\DeclareMathOperator{\E}{\mathbb{E}}

\renewcommand{\P}{\mathbb{P}}

\newcommand{\calF}{\mathcal{F}}

\newcommand{\calL}{\mathcal{L}}
\newcommand{\calU}{\mathcal{U}}
\newcommand{\calT}{\mathcal{T}}
\newcommand{\calG}{\mathcal{G}}
\newcommand{\calA}{\mathcal{A}}
\newcommand{\calB}{\mathcal{B}}

\renewcommand{\bar}[1]{\overline{#1}}

\newcommand{\egaldistr}{{\overset{(d)}{=}}}

\title{Maximal displacement in a branching random walk through~interfaces}
\author{Bastien Mallein\footnote{LPMA, Univ. P. et M. Curie (Paris 6). Research partially supported by the ANR project MEMEMO.}\footnote{DMA, École Normale Supérieure (Paris).}}
\date{\today}

\newcommand{\T}{\mathbf{T}}
\renewcommand{\tilde}[1]{\widetilde{#1}}
\renewcommand{\hat}[1]{\widehat{#1}}

\newcommand{\triangleBRW}[2]{\shade [top color=black!80, bottom color=white] (#1-0.75,-8) -- (#1,-#2-#2) -- (#1+0.75,-8) -- cycle;
	\draw (#1,-#2-#2) node {$\bullet$} ;
	\draw (#1,-7.5) node {$\P_{\cdot,#2}$} ;
}

\newcommand{\bbv}{\mathbf{v}}
\newcommand{\bba}{\mathbf{a}}
\newcommand{\bbb}{\mathbf{b}}
\newcommand{\bbtheta}{\boldsymbol{\theta}}

\begin{document}

\maketitle

\begin{abstract}
In this article, we study a branching random walk in an environment which depends on the time. This time-inhomogeneous environment consists of a sequence of macroscopic time intervals, in each of which the law of reproduction remains constant. We prove that the asymptotic behaviour of the maximal displacement in this process consists of a first ballistic order, given by the solution of an optimization problem under constraints, a negative logarithmic correction, plus stochastically bounded fluctuations.
\end{abstract}

\section{Introduction}
\label{sec:introduction}

The theory of branching processes grew from the seminal work of Galton and Watson to model the dynamic of family names. The Galton-Watson branching process corresponds to a population in which each individual in the generation $n$ independently produces a random number of children, with the same distribution. This process is constructed by recurrence as follows: given $(\xi_{n,k}, n \in \N, k \in \N)$ an i.i.d. array of integer-valued random variables, we write
\[
  Z_0 = 1 \quad \mathrm{and} \quad \forall n \in \N,  Z_{n} = \sum_{j=1}^{Z_{n-1}} \xi_{n,j}.
\]
For $(n,j) \in \N^2$, $\xi_{n,j}$ is the number of children of the $j^\mathrm{th}$ individual alive at generation $n-1$.

A natural development of this model consists of mapping every individual in this Galton-Watson process with a position on the real line. The initial ancestor --i.e. the one individual alive at time $0$-- is positioned at the origin, and the relative position of one individual with respect to its parent is sampled according to an i.i.d. random variable on $\R$. This process is called branching random walk. In greater generality, the displacement of a child does not have to be independent of the displacement of its siblings, or of the number of siblings it has. In this case, the relative position of the children of an individual with respect to their parent forms a point process on $\R$, which characterizes the reproduction of the individual.

In this article, we take interest in time-inhomogeneous branching random walks, in which the reproduction law of individuals depends on the time. Such a time-inhomogeneous branching random walk on $\R$ is a process which starts with one individual located at the origin at time 0, and evolves as follows: at each time $k \in \N$, every individual currently in the process dies, giving birth to a certain number of children, which are positioned around their parent according to independent versions of a point process, whose law may depend on the generation of the parent.

When the law of the point process does not depend on the generation of the individual, and satisfies some integrability conditions, the asymptotic of the maximal displacement is fully known. In the '70s, Hammersley \cite{Ham74}, Kingman \cite{Kin75} and Biggins \cite{Big76} proved this maximal value grows at linear speed almost surely. Hu and Shi \cite{HuS09} exhibited a logarithmic correction in probability, with almost sure fluctuations; while Addario-Berry and Reed \cite{ABR09} showed the tightness of the maximal displacement, shifted around its median. More recently, Aidékon \cite{Aid13} proved the fluctuations converge in law to a random shift of a Gumbel variable.

Fang and Zeitouni \cite{FaZ12a} introduced a time-inhomogeneous branching random walks of length $n\in\N$, defined as follows. At each step, individuals split independently into two children, which move around their parent according to independent Gaussian random variables. During the first $\frac{n}{2}$ units of time, the Gaussian random variables have variance $\sigma_1^2$, while they have variance $\sigma_2^2$ after time $\frac{n}{2}$. The behaviour of this process depends on the sign of $\sigma_2^2-\sigma_1^2$. The asymptotic of the maximal displacement is once again composed by a first ballistic order, a second logarithmic term and fluctuations of order 1; but the logarithmic correction term exhibits a phase transition as $\sigma^2_2$ grows bigger than $\sigma^2_1$.

This result can be extended to more general time-inhomogeneous environments. In this article, we do not assume the displacement of the children to be Gaussian, or independent of the displacement of its siblings. Moreover, we can assume the reproduction law to change more than once in the process. Let $P>0$ be an integer, $0=\alpha_0<\alpha_1< \cdots < \alpha_P=1$ be a partition of $[0,1]$ and $(\calL_p, 1 \leq p \leq P)$ be a family of laws of point processes. We study a time-inhomogeneous branching random walk in which the reproduction law of individuals is equal to $\calL_p$ between time $n \alpha_{p-1}$ and $n\alpha_p$. More precisely, given $n \in \N$ to be the length of the process, we consider a process starting from one individual alive at time $0$ at position $0$; such that at each individual alive at generation $k \in [n \alpha_{p-1}, n\alpha_p)$ reproduces according to an independent point process with law $\calL_p$. We call this process branching random walk through a series of interfaces, as the way individuals reproduce sharply changes at some given times. We prove that in this process, the asymptotic of the maximal displacement is again a first ballistic order plus logarithmic corrections and fluctuations of order 1, under suitable integrability conditions. The value of logarithmic correction is influenced by the path followed by the individual that reaches the maximal position at time $n$.

In this article, $c,C$ are two positive constants, respectively small enough and large enough, which may change from line to line, and depend only on the law of the random variables we consider. For a given sequence of random variables $(X_n, n \geq 1)$, we write $X_n = O_\P(1)$ if the sequence is tensed, i.e. $\lim_{K \to +\infty} \sup_{n \geq 1} \P(|X_n| \geq K) = 0$. Moreover, we always assume the convention $\max \emptyset = - \infty$ and $\min \emptyset = +\infty$, and for $u \in \R$, we write $u_+ = \max(u,0)$, and $\log_+(u) = (\log u)_+$. Finally, $\calC_b$ is the set of continuous bounded functions on $\R$.

In the rest of the introduction, we introduce in Section \ref{subsec:notation} some additional notation on trees, point processes and branching random walks, to give a formal definition of our model in Section \ref{subsec:brwis}. In Section \ref{subsec:assuptions}, we detail the heuristic that can be used to conjecture the value of the first two orders of the asymptotic of $M_n$, before stating our main result in Section \ref{subsec:mainresult}. The rest of the article is devoted to the proof of this result, using the spinal decomposition of the branching random walk, bounds on the probability for a --time-inhomogeneous-- random walk to make an excursion, and some Lagrange multipliers analysis.

\subsection{Definition of the model and notation}
\label{subsec:notation}

\subsubsection{Plane rooted marked trees}

Following the Ulam-Harris notations for trees, we write
\[
  \calU^* = \bigcup_{n \in \N} \N^n \quad \mathrm{and} \quad \calU = \calU^* \cup \{ \emptyset \}
\]
the set of finite sequences of integers, with the convention $\N^0 = \{\emptyset\}$, where $\emptyset$ is the sequence of length $0$, which encodes the root of the trees we consider.

Let $u=(u(1),\ldots u(n)) \in \calU^*$, then $u$ represents the $u(n)^\mathrm{th}$ child of the $u(n-1)^\mathrm{th}$ child of ... of the $u(1)^\mathrm{th}$ child of the initial individual $\emptyset$. We write $|u|=n$ the generation to which $u$ belongs, with convention $|\emptyset|=0$. For any $1 \leq k \leq n$, we set $u_k=(u(1),\ldots u(k))$, and $u_0=\emptyset$. We define the application
\[
  \pi : \begin{array}{rcl}
  \calU^* & \longrightarrow & \calU\\
  (u(1),\ldots u(n)) & \longmapsto & u_{|u|-1} = (u(1),\ldots u(n-1))
  \end{array}
\]
which associates to a vertex $u$ its parent $\pi u$. Note that $\emptyset$ is the only vertex with no parent. For $u, v \in \calU$, we write $u<v$ if there exists $k<|v|$ such that $u=v_k$, or in other words, if $u$ is an ancestor of $v$.

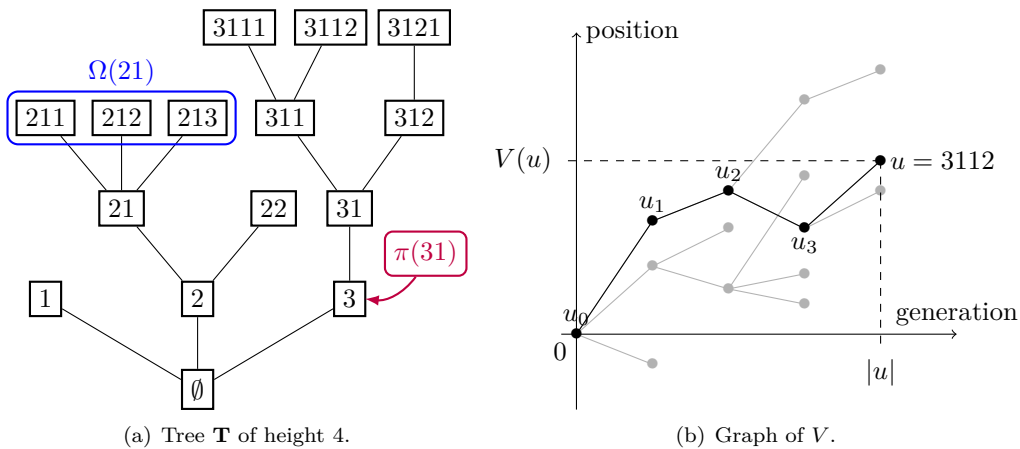
\begin{figure}[ht]
\centering
\subfigure[Tree $\T$ of height 4.]{
\begin{tikzpicture}
  \node [individu] {$\emptyset$} [grow'=up, level distance = 1.2cm, sibling distance = 2cm]
    child { node [individu] {$1$} }
    child { node [individu] {$2$} [sibling distance = 2cm]
      child { node [individu] {$21$} [sibling distance = 1cm]
        child { node [individu] {$211$} }
        child { node [individu] {$212$} }
        child { node [individu] {$213$} }
      }
      child { node [individu] {$22$} }
    }
    child { node [individu] {$3$}
      child { node [individu] {$31$} [sibling distance = 1.7cm]
        child { node [individu] {$311$} [sibling distance = 1.2cm]
          child { node [individu] {$3111$} }
          child { node [individu] {$3112$} }
        }
        child {node [individu] {$312$}
          child { node [individu] {$3121$} }
        }
      }
    };
    \node[draw,rectangle,rounded corners=3pt] [thick, color=purple] (P) at (3,1.8) {$\pi(31)$};
  \draw[->,>=latex, thick, color=purple] (P) to[bend left] (2.2,1.2);
  \draw[thick, color=blue, rounded corners] (-2.5,3.95) rectangle (0.5,3.25) ;
  \draw[thick, color=blue] (-1,3.9) node[above] {$\Omega(21)$} ;
\end{tikzpicture}
}%
\subfigure[Graph of $V$.]{
\begin{tikzpicture}
  \draw[->] (-0.2,0) -- (5,0);
  \draw (5,0) node[above] {generation};
  \draw [->] (0,-1) -- (0,4);
  \draw (0,4) node[right] {position};
  \draw (0,0) node[below left] {$0$};
  \draw[color=black!30] (1,-0.4) -- (0,0) -- (1,0.9) ;
  \draw[color=black!30] (0,0) -- (1,1.5) ;
  \draw[color=black!30] (1,0.9) node{$\bullet$} ;
  \draw[color=black!30] (1,1.5) node{$\bullet$} ;
  \draw[color=black!30] (1,-0.4) node{$\bullet$} ;
  \draw[color=black!30] (2, 0.6) -- (1,0.9) -- (2,1.4) ;
  \draw[color=black!30] (2,0.6) node{$\bullet$};
  \draw[color=black!30] (2,1.4) node{$\bullet$};
  \draw[color=black!30] (1,1.5) -- (2,1.9) ;
  \draw[color=black!30] (2,1.9) node{$\bullet$};
  \draw[color=black!30] (3,0.8) -- (2,0.6) -- (3,0.4) ;
  \draw[color=black!30] (2,0.6) -- (3,2.1) ;
  \draw[color=black!30] (3,0.4) node{$\bullet$};
  \draw[color=black!30] (3,0.8) node{$\bullet$};
  \draw[color=black!30] (3,2.1) node{$\bullet$};
  \draw[color=black!30] (3,1.4) -- (2,1.9) -- (3,3.1) ;
  \draw[color=black!30] (3,1.4) node{$\bullet$};
  \draw[color=black!30] (3,3.1) node{$\bullet$};  
  \draw[color=black!30] (4,1.9)--(3,1.4)--(4,2.3) ;
  \draw[color=black!30] (4,1.9) node{$\bullet$};
  \draw[color=black!30] (4,2.3) node{$\bullet$};  
  \draw[color=black!30] (3,3.1) -- (4,3.5) ;
  \draw[color=black!30] (4,3.5) node{$\bullet$};
  \draw (0,0) -- (1,1.5) -- (2,1.9) -- (3,1.4) -- (4,2.3) ;
  \draw (4,2.3) node {$\bullet$} ;
  \draw (3,1.4) node {$\bullet$} ;
  \draw (2,1.9) node {$\bullet$} ;
  \draw (1,1.5) node {$\bullet$} ;
  \draw (4,2.3) node[right] {$u = 3112$} ;
  \draw [dashed] (4,2.3) -- (-0.2,2.3) ;
  \draw [dashed] (4,2.3) -- (4,-0.2) ;
  \draw (4,-0.2) node[below] {$|u|$} ;
  \draw (-0.2,2.3) node[left] {$V(u)$} ;
  \draw (0,0) node[above] {$u_0$} ;
  \draw (1,1.5) node[above] {$u_1$} ;
  \draw (2,1.9) node[above] {$u_2$} ;
  \draw (3,1.4) node[below] {$u_3$} ;
  \draw[color=black] (0,0) node{$\bullet$};
\end{tikzpicture}
}
\caption{A plane rooted marked tree $(\T,V)$}
\end{figure}

A plane rooted tree $\T$ is a subset of $\calU$ which satisfies the three following properties:
\begin{description}
  \item[(T1)] the root $\emptyset \in \T$ ;
  \item[(T2)] if $u \in \T$ and $u \neq \emptyset$ then $\pi u \in \T$ ;
  \item[(T3)] if $(u(1),\ldots u(n)) \in \T$ and $v \leq u(n)$, then $(u(1),\ldots u(n-1), v) \in \T$.
\end{description}
For example, a Galton-Watson tree can be constructed as follows: given a family $(\xi_{u}, u \in \calU)$ of i.i.d. random variables, we define
\[
  \T = \left\{ u \in \calU : \forall k < |u|, u(k) \leq \xi(u_{k-1}) \right\},
\]
which is indeed a tree. Observe that in this settings, if $u \in \T$ then $\xi(u)$ is the number of children of $u$.

We call height of $\T$ the quantity $\max_{u \in \T} |u|$. All the trees we consider in this article are of finite height. The set $\{ u \in \T : |u|=n\}$ is referred to as the $n^\text{th}$ generation of $\T$, often abbreviated as $\{|u|=n\}$. For a given $u \in \T$, we write $\Omega(u) = \left\{ v \in \T : \pi v = u \right\}$ the set of children of $u$.

A plane rooted marked tree is a pair $(\T,V)$, where $\T$ is a plane rooted marked tree and $V : \T \to \R$. In the context of branching random walks, we refer to $V(u)$ as to the position of individual $u$. The set of plane rooted marked trees is written $\calT$.

\subsubsection{Point processes}

A point process $L$ is a random variable taking values in the set of finite or infinite sequences of real numbers. Once again, the empty sequence is written $\emptyset$. The point processes we consider in this article admit a maximum and have no accumulation point. Therefore, we write $L=(\ell_1,\ldots, \ell_N)$, where $\ell_1\geq \ell_2 \geq \cdots$ is the set of points in $L$, with the convention $\ell_{+\infty}=-\infty$ and $N$ is a random variable taking values in $\Z_+ \cup \{+\infty\}$, which represents the total number of points in $L$. We write $\calL$ the law of $L$. Using the same vocabulary as in Galton-Watson processes, we say that $\calL$ never gets extinct and has supercritical offspring if
\begin{equation}
  \label{eqn:breeding}
  \P(L= \emptyset) = \P(N=0) = 0 \quad \mathrm{and} \quad \E\left( \sum_{\ell \in L} 1 \right) = \E(N) > 1.
\end{equation}

For any $\theta \geq 0$, we write $\kappa(\theta) = \log \E\left[ \sum_{\ell \in L} e^{\theta \ell} \right]$ the log-Laplace transform of $\calL$, and for all $a \in \R$, $\kappa^*(a) = \sup_{\theta > 0} \left[ \theta a - \kappa(\theta) \right]$ its Fenchel-Legendre transform. Let $f : \R_+ \to \R \cup \{+\infty\}$ be a convex function, and $f^*$ its transform. If $f^*$ is differentiable at point $x$ then
\begin{equation}
  \label{eqn:legendreestimate}
  f^*(x) = (f^*)'(x) x - f\left( (f^*)'(x) \right).
\end{equation}

\subsubsection{Branching random walk in time-inhomogeneous environment}

A branching random walk is a random variable taking values in $\calT$ the set of rooted marked trees. Let $n \in \N$ and $(\calL_1,\ldots, \calL_n)$ be a family of point processes laws, which we call the environment of the branching random walk. The law of the time-inhomogeneous branching random walk $(\T,V)$ of length $n$ with environment $(\calL_1,\ldots \calL_n)$ is characterized by the three following properties
\begin{description}
  \item[(BRWtie1)] $V(\emptyset)=0$;
  \item[(BRWtie2)] $\left\{ \left( V(v)-V(u), v \in \Omega(u) \right) u \in \T \right\}$ is a family of independent point processes;
  \item[(BRWtie3)] $\left( V(v) - V(u), v \in \Omega(u)\right)$ has law $\calL_{|u|+1}$, where $\calL_{n+1} = \delta_{\emptyset}$.
\end{description}
This branching random walk can be constructed as follows. We consider a family of independent point processes $\left\{ L^u, u \in \calU, |u| \leq n-1 \right\}$, where $L^u$ has law $\calL_{|u|+1}$. For any $u \in \calU$ with $|u|<n$, we write $L^u = (\ell^u_1,\ldots \ell^u_{N(u)})$. The plane rooted tree which represents the genealogy of the population is
\[
  \T = \left\{ u \in \calU : |u| \leq n, \forall k \leq |u|-1, u(k+1) \leq N(u_k)\right\}.
\]
We observe that $\T$ is a --time-inhomogeneous-- Galton-Watson tree, with reproduction law at generation $k$ given by the number of points in a point process of law $\calL_k$. We set $V(\emptyset) = 0$ and, for $u \in \T$ with $|u|=k$,
\begin{equation*}
  V(u) := V(\pi u) + \ell^{\pi u}_{u(k)} = \sum_{j=0}^{k-1} \ell^{u_j}_{u(j+1)}.
\end{equation*}
For $u \in \T$, we often call path or trajectory of $u$ the sequence $(V(u_0), V(u_1), \ldots V(u))$ of positions of the ancestors of $u$. Finally, we write $M_n = \max_{|u|=n} V(u)$ the maximal displacement in the branching random walk at generation $n$.

\subsubsection{Branching random walk through a series of interfaces}
\label{subsec:brwis}

In this article, we take interest in branching random walks through interfaces. In this model, the time-inhomogeneous environment consists of a series of macroscopic stages. We set $P \in \N$ the number of such stages, $0=\alpha_0<\alpha_1 < \cdots < \alpha_P=1$ the times at which the interfaces occur, and $(\calL_p, p \leq P)$ a $P$-uple of laws of point processes.

For $n \in \N$ and $p \leq P$, we write $\alpha^{(n)}_p = \floor{n \alpha_p}$. The branching random walk through a series of interfaces --BRWis for short-- of length $n$ is a branching random walk in time-inhomogeneous environment, in which individuals alive at generation $k$ reproduce according to the law $\calL_p$ for all $\alpha^{(n)}_{p-1} \leq k < \alpha^{(n)}_p$. We write $(\T^{(n)}, V^{(n)})$ such a branching random walk. When the value of $n$ is clear in the context, we often omit the superscripts to make the notations lighter.

The law of $(\T^{(n)},V^{(n)})$ is characterized by the three following properties
\begin{description}
  \item[(BRWis1)] $V^{(n)}(\emptyset)=0$ ;
  \item[(BRWis2)] $\left\{ \left( V^{(n)}(v)-V^{(n)}(u), v \in \Omega(u) \right) u \in \T^{(n)} \right\}$ is a family of independent point processes ;
  \item[(BRWis3)] $\left( V^{(n)}(v) - V^{(n)}(u), v \in \Omega(u)\right)$ has law $\calL_{p}$ if $n\alpha_{p-1} \leq |u| < n\alpha_p$ and is empty otherwise.
\end{description}

\begin{remark}
By splitting the first time-interval of the BRWis into three pieces, we always assume that the number $P$ of stages we consider is greater than or equal to $3$ in the rest of the article. In particular, the results we obtain here hold for time-homogeneous branching random walks.
\end{remark}

\subsection{Assumptions and main result}
\label{subsec:assuptions}

We fix an integer $P$, a sequence $0=\alpha_0<\alpha_1 < \cdots < \alpha_P=1$ and a family $(\calL_p, p \leq P)$ of points processes laws. We write $\kappa_p$ the log-Laplace transform of $\calL_p$, and $\kappa^*_p$ its Fenchel-Legendre transform. We introduce some well-known branching random walk estimates, and use them to build heuristics for the comportment of the BRWis, before stating the main result of this article.

\subsubsection{Some well-known estimates for a time-homogeneous branching random walk}

We list some classical branching random walk results, that can be found in \cite{Big10}. Let $p \leq P$, we consider a time-homogeneous branching random walk $(\T_p,V_p)$, in which individuals reproduce according to law $\calL_p$. We write $M_{p,n} = \max_{|u|=n} V_p(u)$ its maximal displacement at time $n$. If there exists $\theta > 0$ such that $\kappa_p(\theta)<+\infty$, we set
\begin{equation}
  \label{eqn:speeddef}
  v_p = \inf_{\theta > 0} \frac{\kappa_p(\theta)}{\theta} = \sup\{a \in \R : \kappa^*(a) \leq 0 \}.
\end{equation}
As $\lim_{n \to +\infty} \frac{M_{p,n}}{n} = v_p$ a.s, $v_p$ is called the {\em speed} of the branching random walk. Under the assumption
\begin{equation}
  \label{eqn:regularity}
  \forall p \leq P, \exists \bar{\theta}_p \in \R_+ : \bar{\theta}_p \kappa'_p(\bar{\theta}_p) - \kappa_p(\bar{\theta}_p)=0,
\end{equation}
we have $v_p = \kappa'_p(\bar{\theta}_p)$. Moreover, the function $\kappa^*_p$ is linked to the density of individuals present in the $n^\mathrm{th}$ generation. As proved in \cite{Big77}, we have
\begin{equation}
  \label{eqn:density}
  \begin{cases}
    \forall a < v_p, \lim_{n \to +\infty} \frac{1}{n} \log \sum_{|u|=n} \ind{V_p(u) \geq na} = -\kappa_p^*(a) \quad \mathrm{a.s.}\\
    \forall a > v_p, \lim_{n \to +\infty} \frac{1}{n} \log \P\left[ \exists |u| = n: V_p(u) \geq na \right] = - \kappa_p^*(a).
  \end{cases}
\end{equation}
With high probability, there is no individual above $v_p$, and there is an exponentially large number of individuals above $n(v_p-\epsilon)$. More precisely, by equation \eqref{eqn:density}, $e^{-n \kappa^*(a)}$ is either an approximation of the number of individuals alive at time $n$ in a neighbourhood of $na$, or of the probability to observe at least one individual around $na$ at time $n$, depending on the sign of $\kappa^*(a)$.

\subsubsection{Heuristics for the maximal displacement}

We now consider the BRWis $(\T,V)$. Given $\bba = (a_p, p \leq P) \in \R^P$ --in the rest of the article, we write in bold letters real $P$-uples-- we take interest in the number of individuals alive at time $n$ such that for any $p < P$, their ancestor at time $\alpha^{(n)}_p$ were close to $n\sum_{k=1}^p a_k (\alpha^{(n)}_{k}-\alpha^{(n)}_{k-1})$. For every such individual, we say that it ``follows the path driven by $\bba$''.

Using \eqref{eqn:density}, we know there are $e^{-\alpha^{(n)}_1 \kappa^*_1(a_1)}$ individuals alive at time $\alpha^{(n)}_1$ around $\alpha^{(n)}_1 a_1$ if $\kappa^*_1(a_1)<0$, and none otherwise. Each individual starts an independent branching random walk from $\alpha^{(n)}_1 a_1$. Applying the law of large numbers, we expect $e^{-\alpha^{(n)}_1 \kappa^*_1(a_1) - (\alpha^{(n)}_2-\alpha^{(n)}_1) \kappa^*_2(a_2)}$ descendants at time $\alpha^{(n)}_2$ at position $\alpha^{(n)}_1 a_1 + (\alpha^{(n)}_2 - \alpha^{(n)}_1) a_2$. More generally, we write
\[
  K^* : \begin{array}{rcl}
  \R^P & \to & \R^P\\
  \bba & \mapsto & \left( \sum_{q=1}^p (\alpha_q - \alpha_{q-1}) \kappa^*_1(a_q), p \leq P\right).
  \end{array}
\]
For any $p \leq P$, we expect $e^{-n K^*(a)_p}$ individuals who followed the path driven by $\bba$ until time $\alpha^{(n)}_p$.

Let $\bba \in \R^P$, if for all $p \leq P$, $K^ *(\bba)_p \leq 0$, we expect at time $n$ about $e^{- n K^*(\bba)_P}$ individuals who followed the path driven by $a$. In particular, this means that there is at least one individual above $\sum_{p=1}^P (\alpha_p-\alpha_{p-1})a_p$. On the other hand, if there exists $p_0 \leq P$ such that $K^*(\bba)_{p_0} > 0$, then with high probability, no individual alive at time $\alpha^{(n)}_{p_0}$ followed this path.

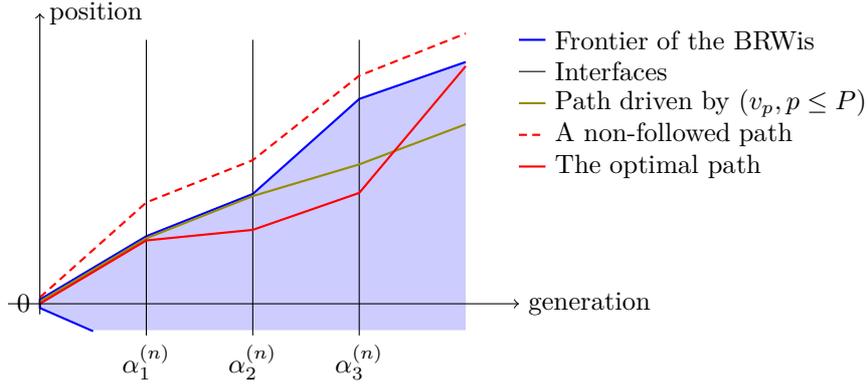
\begin{figure}[ht]
\centering
\begin{tikzpicture}[scale=0.7]
  \fill [color=blue!20] (8,-0.5) -- (1,-0.5)--(0,-0.08)--(0,0.08) -- (2,1.28) -- (4,2.08) -- (6,3.88) -- (8,4.58) -- cycle ;
  \draw [thick, color=blue] (1,-0.52)--(0,-0.08)--(0,0.08) -- (2,1.28) -- (4,2.08) -- (6,3.88) -- (8,4.58) ;
  \draw (2,-0.6) -- (2,5) ;
  \draw (2,-0.6) node[below]{$\alpha^{(n)}_1$} ;
  \draw (4,-0.6) -- (4,5) ;
  \draw (4,-0.6) node[below]{$\alpha^{(n)}_2$} ;
  \draw (6,-0.6) -- (6,5) ;
  \draw (6,-0.6) node[below]{$\alpha^{(n)}_3$} ;
  \draw [thick, color=olive] (0,0.04) -- (2,1.24) -- (4,2.04) -- (6,2.64) -- (8,3.4) ;
  \draw [thick,densely dashed, color=red] (0,0.12) -- (2,1.92) -- (4,2.72) -- (6,4.32) -- (8,5.12) ;
  \draw [thick, color=red] (0,0) -- (2,1.2) -- (4,1.4) -- (6,2.1) -- (8,4.5) ;
  
  \draw [thick, color=blue] (9,5) -- (9.5,5) ;
  \draw (9.5,5) node[right] {Frontier of the BRWis};
  \draw (9,4.4) -- (9.5,4.4) ;
  \draw (9.5,4.4) node[right] {Interfaces};
  \draw [thick, color=olive] (9,3.8) -- (9.5,3.8) ;
  \draw (9.5,3.8) node[right] {Path driven by $(v_p, p \leq P)$};
  \draw [thick, densely dashed, color=red] (9,3.2) -- (9.5,3.2) ;
  \draw (9.5,3.2) node[right] {A non-followed path};
  \draw [thick, color=red] (9,2.6) -- (9.5,2.6) ;
  \draw (9.5,2.6) node[right] {The optimal path};
  \draw[->] (-0.6,0) -- (9,0);
  \draw (9,0) node[right] {generation};
  \draw [->] (0,-0.2) -- (0,5.5);
  \draw (0,5.5) node[right] {position};
  \draw (0,0) node[left] {$0$};
\end{tikzpicture}
\caption{Different path of interest in the BRWis.}
\end{figure}

We write $\calR = \left\{ \bba \in \R^P : \forall p \leq P, K^*(\bba)_p \leq 0 \right\}$. Following the heuristic, we expect to find individuals alive in the process at time $n$ around position $nu$ if and only if $u = \sum (\alpha_p-\alpha_{p-1}) a_p$ for some $\bba \in \calR$. We set
\begin{equation}
  \label{eqn:speedisdef}
  v_\mathrm{is} = \sup_{\bba \in \calR}  \sum_{p=1}^P (\alpha_p-\alpha_{p-1}) a_p
\end{equation}
which we prove to be the speed of the BRWis.

\subsubsection{The optimization problem}

According to this heuristic, if the BRWis verifies
\begin{equation}
  \label{eqn:optimization_problem}
  \exists \bba \in \calR : v_\mathrm{is} = \sum_{p=1}^P (\alpha_p - \alpha_{p-1}) a_p,
\end{equation}
then the path followed by the rightmost individual until time $n$ is driven by the optimal solution $\bba$. Under the additional assumption
\begin{equation}
\label{eqn:differentiable}
  \forall p \leq P, \forall a \in \R,  \kappa^*_p \text{ is differentiable at point } a \text{ or } \kappa^*_p(a)=+\infty,
\end{equation}
this optimal solution satisfies some interesting properties. To guarantee existence and/or uniqueness of the solutions of \eqref{eqn:optimization_problem}, we need to introduce additional integrability assumptions, such as
\begin{equation}
  \label{eqn:finite_reproduction}
  \forall p \leq P, \kappa_p(0) \in (0, +\infty) \quad \mathrm{and} \quad \kappa'_p(0) \text{ exists}.
\end{equation}

\begin{proposition}
\label{prop:optimization_problem}
If point processes $\calL_1,\ldots \calL_P$ verify \eqref{eqn:breeding}, under assumption \eqref{eqn:differentiable}, $\bba \in \calR$ is a solution of \eqref{eqn:optimization_problem} if and only if, writing $\theta_p = \left( \kappa^*_p \right)'(a_p)$, we have
\begin{enumerate}
  \item $\bbtheta$ is non-decreasing and positive ;
  \item if $K^*(\bba)_p<0$, then $\theta_{p+1}=\theta_p$ ;
  \item $K^*(\bba)_P=0$.
\end{enumerate}

Under the conditions \eqref{eqn:regularity} and \eqref{eqn:differentiable}, there exists at most one solution to \eqref{eqn:optimization_problem}.

Under the conditions \eqref{eqn:differentiable} and \eqref{eqn:finite_reproduction}, there exists at least one solution to \eqref{eqn:optimization_problem}.
\end{proposition}

The proof of this result, which is a direct application of the theory of Lagrange multipliers, is postponed to Appendix \ref{app:optimization_problem}. Despite the fact that this would be a natural candidate, the path driven by $\bbv := (v_1,\ldots, v_P)$ is not always the optimal solution. For example, if there exists $p \leq P-1$ such that $\bar{\theta}_p > \bar{\theta}_{p+1}$, Proposition \ref{prop:optimization_problem} proves that $\bbv$ is not the solution. Loosely speaking, in this case, the path of the rightmost individual at time $n$ does not stay close to the boundary of the branching random walk at all time.

\begin{remark}
\label{rem:decreasing_variance}
On the other hand, under assumptions \eqref{eqn:regularity} and \eqref{eqn:differentiable}, if $\mathbf{\bar{\theta}}$ is positive and non-decreasing, then $\bbv$ is indeed the optimal solution. In this case, $\bbv$ satisfies the first assumption of Proposition \ref{prop:optimization_problem}, and the two others are an easy consequence of $K^*(\bbv)_p=0$ for any $p \leq P$. This situation corresponds, in Gaussian settings, to branching random walks with decreasing variance. In this situation, the rightmost individual at time $n$ stays at any time $k<n$ within range $O(n^{1/2})$ from the frontier of the BRWis.
\end{remark}

\subsubsection{On the logarithmic correction}

We discuss the heuristic for the logarithmic correction of the BRWis. For a time-homogeneous branching random walk with reproduction law $\calL_p$, under assumption \eqref{eqn:regularity} and some additional integrability estimates, we have
\[
  M_n^{(p)} = nv_p - \frac{3}{2 \bar{\theta}_p} \log n +O_\P(1),
\]
and the second order can be directly related, as it is underlined in \cite{AiS10}, to the following estimate for a random walk with finite variance,
\[
  \log \P\left[ S_n \leq \E(S_n) + 1, S_j \geq \E(S_j), j \leq n \right] \sim_{n \to +\infty} -\frac{3}{2}\log n.
\]
In effect, the path followed by the rightmost individual at time $n$ made an excursion below the frontier of the branching random walk.

A similar condition holds for BRWis, the path leading to the rightmost individual at time $n$ stays below the frontier of the branching random walk at any time $k \leq n$. Note that if $K^*(\bba)_p=0$, then the optimal path is at distance $o(n)$ from the frontier of the branching random walk. Moreover, for any $p$ such that $\theta_{p+1}>\theta_p$, we prove the ancestor at time $\alpha^{(n)}_p$ of the rightmost individual at time $n$ was within distance $O(1)$ of the frontier. The logarithmic correction is a sum of terms related to the difficulty for a random walk to stay below the boundary of the branching random walk, and hit at time $n$ this boundary.

From now on, $\bba$ stands for the optimal solution of \eqref{eqn:speedisdef}, and $\theta_p = (\kappa^*_p)'(a_p)$. We denote the number of different values taken by $\bbtheta$ by $T = \# \{ \theta_p, p \leq P \}$ and set $\phi_1 < \phi_2 < \cdots < \phi_T$ the distinct values taken by $\bbtheta$, listed in the increasing order. For any $t \leq T$, we set $f_t = \min\{p \leq P : \theta_p = \phi_t\}$ and $l_t = \max\{p \leq P : \theta_p = \phi_t\}$. Observe that for any $p \in [f_t,l_t]$, we have $\theta_p=\phi_t$. We write
\begin{equation}
  \label{eqn:logisdef}
  \lambda = \sum_{t=1}^T \frac{1}{2\phi_t} \left[ \ind{K^*(\bba)_{f_t} = 0} + 1 + \ind{K^*(\bba)_{l_t-1}=0} \right]
\end{equation}
with the convention $K^*(\bba)_0=0$. Condition $K^*(\bba)_{f_t}=0$ means that between times $\alpha^{(n)}_{f_t-1}$ and $\alpha^{(n)}_{f_t}$, the optimal path stays close to the frontier of the BRWis, which has a cost of order $\frac{1}{2} \log n$ by the ballot theorem (see Section \ref{sec:randomwalk}). Moreover, each time the value of $\bbtheta$ changes, the optimal path is localized in a window of width $O(1)$, which has cost $\frac{1}{2}\log n$ by the local limit theorem. We prove that under some good integrability conditions $M_n \approx nv-\lambda \log n$. We observe that $\lambda \geq \frac{1}{2\phi_1}>0$. If $P=T=1$, then $\lambda = \frac{3}{2\phi_1}$, which is consistent with the results of Hu--Shi and Addario-Berry--Reed.

\subsubsection{The asymptotic of the maximal displacement in the BRWis}
\label{subsec:mainresult}

We recall that $\bba$ is the solution of \eqref{eqn:optimization_problem}. We write
\begin{equation}
  \label{eqn:strongfrontier}
  B=\left\{ p \leq P : K^*(\bba)_{p-1}=K^*(\bba)_p=0\right\},
\end{equation}
such that for any $k \in \cup_{p \in B} [\alpha^{(n)}_{p-1}, \alpha^{(n)}_p]$, the path leading to the the rightmost individual is within distance $o(n)$ from the frontier of the branching random walk. For any $p \leq P$, we introduce the random variable 
\begin{equation}
  \label{eqn:definexp}
  X_p = \sum_{\ell \in L_p} e^{\theta_p \ell} \quad \mathrm{and} \quad \tilde{X}_p = \sum_{\ell \in L_p} \ell e^{\theta_p \ell},
\end{equation}
and we assume the following integrability conditions for the point processes:
\begin{equation}
  \label{eqn:variance}
  \sup_{p \leq P} \E\left[ \sum_{\ell \in \calL_p} \ell^2 e^{\theta_p \ell}\right]< +\infty,
\end{equation}
\begin{equation}
  \label{eqn:integrability}
  \sup_{p \in B} \E\left[ X_p \left(\log_+ \tilde{X}_p\right)^2 \right] + \sup_{p \in B^c} \E\left[X_p \log_+ X_p\right] < +\infty
\end{equation}

The following theorem is the main result of the article.
\begin{theorem}
\label{thm:main}
If $\calL_1, \cdots \calL_p$ satisfy \eqref{eqn:breeding}, under assumptions \eqref{eqn:optimization_problem}, \eqref{eqn:differentiable}, \eqref{eqn:variance} and \eqref{eqn:integrability}, we have
\[ M_n  =  nv_\mathrm{is} - \lambda \log n + O_\P(1).\]
\end{theorem}

The rest of the article is organized as follows. In Section \ref{sec:spinaldecomposition}, we introduce the spinal decomposition, that links the additive moments of the branching random walk with random walk estimates. In Section \ref{sec:randomwalk}, we compute upper and lower bounds for the probability for a time-inhomogeneous random walk to make an excursion above a given curve. In Section \ref{sec:tailbehaviour}, we give a tight estimate of the tail of $M_n$, which is enough to prove Theorem \ref{thm:main}, using a standard cutting argument. We discuss in Section \ref{sec:twoenvi} consequences of this result for a branching random walk with one interface.

\section{Spinal decomposition of the time-inhomogeneous branching random walk}
\label{sec:spinaldecomposition}

This section is devoted to the proof of a time-inhomogeneous version of the well-known spinal decomposition of the branching random walk. This result consists of two ways of describing a size-biased version of the law of the branching random walk. The spinal decomposition has been introduced to study Galton-Watson processes in \cite{LPP95}. This result is adapted for the first time in \cite{Lyo97} to the branching random walk settings.

\subsection{The size-biased law of the branching random walk}

Let $n \geq 1$ and $(\calL_k, k \leq n)$ be a sequence of point processes laws which forms the environment of a time-inhomogeneous branching random walk $(\T,V)$. For all $x \in \R$ we set $\P_x$ the law on $\calT$ of the marked tree $(\T,V+x)$, and $\E_x$ the corresponding expectation.

We write $\kappa_k(\theta)$ for the log-Laplace transform of $\calL_k$ and we assume there exists $\theta>0$ such that for any $k \leq n$ we have $\kappa_k(\theta)<+\infty$. Let $W_n = \sum_{|u|=n} \exp\left( \theta V(u) - \sum_{j=1}^n \kappa_j(\theta) \right)$. We observe that $W_n > 0, \P_x-\mathrm{a.s.}$ and $\E_x(W_n)=e^x$. We define the law
\begin{equation}
  \label{eqn:size_biaised_law}
  \bar{\P}_x = e^{-\theta x} W_n \cdot \P_x.
\end{equation} 
The spinal decomposition consists of an alternative construction of the law $\bar{\P}_a$, as the projection of a law on the set of planar rooted marked trees with spine, which we define below.

\subsection{A law on plane rooted marked trees with spine}

Let $(\T,V) \in \calT$ be a tree of height $n$, and $w \in \{ u \in \T : |u|=n\}$ an individual alive at the $n^\text{th}$ generation. The triplet $(\T,V, w)$ is a plane rooted marked tree with spine of length $n$. The spine of a tree is a distinguished path of length $n$ linking the root and the $n^\mathrm{th}$ generation. The set of marked trees with spine of height $n$ is written $\hat{\calT}_n$. On this set, we define the three following filtrations,
\begin{align*}
  \forall k \leq n, \hat{\calF}_k &= \sigma\left( u,V(u), u \in \T, |u| \leq k \right) \vee \sigma(w_j, j \leq k) \quad \mathrm{and} \quad \hat{\calF} = \hat{\calF}_n\\
  \forall k \leq n, \calF_k &= \sigma\left( u, V(u) : u \in \T, |u| \leq k\right) \quad \mathrm{and} \quad \calF = \calF_n \\
  \forall k \leq n, \calG_k &= \sigma\left(w_j, V(w_j) : j \leq k \right) \vee \sigma\left(u,V(u), u \in \Omega(w_j), j < k\right) \quad \mathrm{and} \quad \calG = \calG_n.
\end{align*}
The filtration $\calF$ is the information of the marked tree, obtained by forgetting the spine, $\calG$ is the sigma-field of the knowledge of the spine and its children only, and $\hat{\calF}=\calF \vee \calG$ is the natural filtration of the branching random walk with spine.

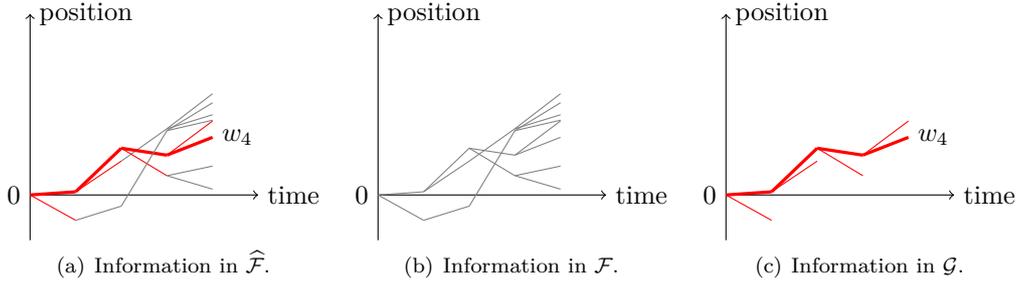
\begin{figure}[ht]
\centering
\subfigure[Information in $\hat{\calF}$.]{
\begin{tikzpicture}[scale=0.6]
  \draw[->] (0,0) -- (5,0);
  \draw (5,0) node[right] {time};
  \draw [->] (0,-1) -- (0,4);
  \draw (0,4) node[right] {position};
  \draw (0,0) node[left] {$0$};
  \draw[very thick, color=red] (0,0) -- (1,0.0696808149084) ;
  \draw[color=red] (0,0) -- (1,-0.5601910912) ;
  \draw[color=red] (1,0.0696808149084) -- (2,0.747674304381) ;
  \draw[very thick, color=red] (1,0.0696808149084) -- (2,1.03299281017) ;
  \draw[color=black!50] (1,-0.5601910912) -- (2,-0.247309772849) ;
  \draw[color=black!50] (2,0.747674304381) -- (3,1.46321068111) ;
  \draw[color=red] (2,1.03299281017) -- (3,0.42544586468) ;
  \draw[very thick, color=red] (2,1.03299281017) -- (3,0.879642841798) ;
  \draw[color=black!50] (2,-0.247309772849) -- (3,1.42114122831) ;
  \draw[color=black!50] (3,1.46321068111) -- (4,2.23774687038) ;
  \draw[color=black!50] (3,1.46321068111) -- (4,1.77253462079) ;
  \draw[color=black!50] (3,0.42544586468) -- (4,0.63777040165) ;
  \draw[color=black!50] (3,0.42544586468) -- (4,0.1278599898) ;
  \draw[color=red] (3,0.879642841798) -- (4,1.6344174291) ;
  \draw[very thick, color=red] (3,0.879642841798) -- (4,1.27365751116) ;
  \draw  (4,1.27365751116) node[right] {$w_4$} ;
  \draw[color=black!50] (3,1.42114122831) -- (4,2.03704867124) ;
  \draw[color=black!50] (3,1.42114122831) -- (4,1.65244242523) ;
\end{tikzpicture}
}%
\subfigure[Information in $\calF$.]
{
\begin{tikzpicture}[scale=0.6]
  \draw[->] (0,0) -- (5,0);
  \draw (5,0) node[right] {time};
  \draw [->] (0,-1) -- (0,4);
  \draw (0,4) node[right] {position};
  \draw (0,0) node[left] {$0$};
  \draw[color=black!50] (0,0) -- (1,0.0696808149084) ;
  \draw[color=black!50] (0,0) -- (1,-0.5601910912) ;
  \draw[color=black!50] (1,0.0696808149084) -- (2,0.747674304381) ;
  \draw[color=black!50] (1,0.0696808149084) -- (2,1.03299281017) ;
  \draw[color=black!50] (1,-0.5601910912) -- (2,-0.247309772849) ;
  \draw[color=black!50] (2,0.747674304381) -- (3,1.46321068111) ;
  \draw[color=black!50] (2,1.03299281017) -- (3,0.42544586468) ;
  \draw[color=black!50] (2,1.03299281017) -- (3,0.879642841798) ;
  \draw[color=black!50] (2,-0.247309772849) -- (3,1.42114122831) ;
  \draw[color=black!50] (3,1.46321068111) -- (4,2.23774687038) ;
  \draw[color=black!50] (3,1.46321068111) -- (4,1.77253462079) ;
  \draw[color=black!50] (3,0.42544586468) -- (4,0.63777040165) ;
  \draw[color=black!50] (3,0.42544586468) -- (4,0.1278599898) ;
  \draw[color=black!50] (3,0.879642841798) -- (4,1.6344174291) ;
  \draw[color=black!50] (3,0.879642841798) -- (4,1.27365751116) ;
  \draw[color=black!50] (3,1.42114122831) -- (4,2.03704867124) ;
  \draw[color=black!50] (3,1.42114122831) -- (4,1.65244242523) ;
\end{tikzpicture}
}%
\subfigure[Information in $\calG$.]
{
\begin{tikzpicture}[scale=0.6]
  \draw[->] (0,0) -- (5,0);
  \draw (5,0) node[right] {time};
  \draw [->] (0,-1) -- (0,4);
  \draw (0,4) node[right] {position};
  \draw (0,0) node[left] {$0$};
  \draw[very thick, color=red] (0,0) -- (1,0.0696808149084) ;
  \draw[color=red] (0,0) -- (1,-0.5601910912) ;
  \draw[color=red] (1,0.0696808149084) -- (2,0.747674304381) ;
  \draw[very thick, color=red] (1,0.0696808149084) -- (2,1.03299281017) ;
  \draw[color=red] (2,1.03299281017) -- (3,0.42544586468) ;
  \draw[very thick, color=red] (2,1.03299281017) -- (3,0.879642841798) ;
  \draw[color=red] (3,0.879642841798) -- (4,1.6344174291) ;
  \draw[very thick, color=red] (3,0.879642841798) -- (4,1.27365751116) ;
  \draw  (4,1.27365751116) node[right] {$w_4$} ;
\end{tikzpicture}
}%
\caption{The graph of a plane rooted marked tree with spine; and the filtrations of $\hat{\calT}$.}
\end{figure}

We now introduce a law $\hat{\P}_x$ on $\hat{\calT}_n$. For any $k \leq n$, we write $\hat{\calL}_k = \left( \sum_{\ell \in L} e^{\theta \ell - \kappa_k(\theta)} \right) \cdot \calL_k$, a law of a point process with Radon-Nikod\'ym derivative with respect to $\calL_k$, and we write $\hat{L}_k = (\hat{\ell}_k(j), j \leq N_k)$ an independent point processes of law $\hat{\calL}_k$. Conditionally on $(\hat{L}_k, k \leq n)$, we choose, for every $k \leq n$, $w(k) \leq N_k$ independently at random, such that
\[
  \P\left(w(k) = h\left|\hat{L}_k, k \leq n\right. \right) = \ind{h \leq N_k} \frac{e^{\theta \ell_k(h)}}{\sum_{j \leq N_k} e^{\theta \ell_k(j)}}.
\]
We denote by $w_n \in \calU$ the sequence $(w(1), \ldots w(n))$.

Let $\left\{L^u, u \in \calU, |u| \leq n\right\}$ be a family of independent point processes such that $L^{w_k} = \hat{L}_{k+1}$, and if $u \neq w_{|u|}$, then $L^u$ has law $\calL_{|u|+1}$. For any $u \in \calU$ such that $|u| \leq n$, we write $L^u = (\ell^u_1,\ldots \ell^u_{N(u)})$. We construct the random tree
\[
  \T = \left\{ u \in \calU : |u| \leq n, \forall 1 \leq k \leq |u|, u(k) \leq N(u_{k-1})\right\},
\]
and function $V : u \in \T \mapsto \sum_{k=1}^{|u|} \ell^{u_{k-1}}_{u(k)}$. For any $x \in \R$, the law of $(\T, x + V, w_n) \in \hat{\calT}_n$ is written $\hat{\P}_x$, and the corresponding expectation is $\hat{\E}_x$.

The marked tree with spine $(\T,x+V,w_n)$ is called branching random walk with spine, and can be constructed as a process in the following manner. It starts with a unique individual positioned at $x$ at time 0, which is the ancestral spine $w_0$. At each time $k < n$, every individual alive at generation $k$ dies. Each of these individuals gives birth to children, which are positioned around their parent according to an independent point process. If the parent is $w_k$, then the law of this point process is $\hat{\calL}_k$, otherwise it is $\calL_k$. Individual $w_{k+1}$ is then chosen at random among the children $u$ of $w_k$, with probability proportional to $e^{\theta V(u)}$. At time $n$, individuals die without children.

In the rest of the article, we write $\P_{x,k}$ for the law of the time-inhomogeneous branching random walk of length $n-k$ starting from $x$ with environment $(\calL_{k+1}, \ldots \calL_n)$. We observe that conditionally on $\calG_k$, the branching random walks of the descendants of the children of $w_k$ are independent, and the branching random walk of the children of $u \in \Omega(w_k)$ has law $\P_{V(u),k+1}$.

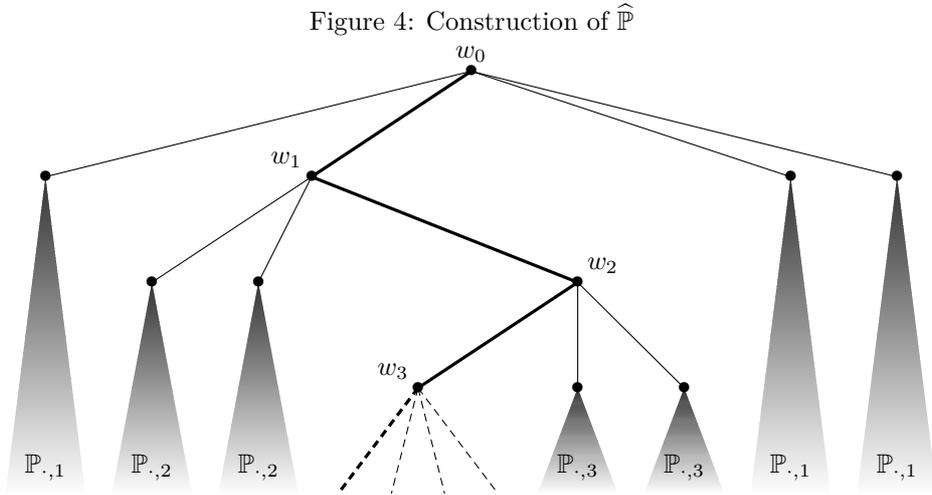
\begin{figure}[ht]
\centering
\caption{Construction of $\hat{\P}$}
\begin{tikzpicture}[scale=0.7]
	\draw [very thick] (0,0) node {$\bullet$} ;
	\draw (0,0) node[above] {$w_0$} ;
	\draw [very thick] (0,0) -- (-3,-2) ;
	\draw (0,0) -- (-8,-2) ;
	\triangleBRW{-8}{1}
	\draw (0,0) -- (6,-2) ;
	\triangleBRW{6}{1}
	\draw (0,0) -- (8,-2) ;
	\triangleBRW{8}{1}
	
	\draw [very thick] (-3,-2) node{$\bullet$} ;
	\draw (-3,-2) node[above left] {$w_1$} ;
	\draw [very thick] (-3,-2) -- (2,-4) ;
	\draw (-3,-2) -- (-6,-4) ;
	\triangleBRW{-6}{2}
	\draw (-3,-2) -- (-4,-4) ;
	\triangleBRW{-4}{2}
	
	\draw [very thick] (2,-4) node{$\bullet$} ;
	\draw (2,-4) node[above right] {$w_2$} ;
	\draw [very thick] (2,-4) -- (-1,-6);
	\draw (2,-4) -- (2,-6);
	\triangleBRW{2}{3}
	\draw (2,-4) -- (4,-6) ;
	\triangleBRW{4}{3} ;
	
	\draw [very thick] (-1,-6) node{$\bullet$} ;
	\draw (-1,-6) node[above left]{$w_3$} ;
	\draw [densely dashed, very thick] (-1,-6) -- (-2.5,-8) ;
	\draw [densely dashed] (-1,-6) -- (-1.5,-8) ;
	\draw [densely dashed] (-1,-6) -- (-0.5,-8) ;
	\draw [densely dashed] (-1,-6) -- (0.5,-8) ;
\end{tikzpicture}
\end{figure}

\subsection{The spinal decomposition}

The following result links the laws $\hat{\P}_x$ and $\bar{\P}_x$ and is the time-inhomogeneous version of the spinal decomposition.
\begin{proposition}[Spinal decomposition]
\label{prop:spinaldecomposition}
For any $x \in \R$, we have
\begin{equation}
  \label{eqn:spinaldecomposition}
  \bar{\P}_x = \left.\hat{\P}_x\right|_\calF.
\end{equation}
Moreover, for any $|u|=n$, we have
\begin{equation}
  \label{eqn:distribspine}
  \hat{\P}_x(w_n = u |\calF) = \frac{\exp\left( \theta V(u) - \sum_{k=1}^n \kappa_k(\theta)\right)}{W_n}.
\end{equation}
\end{proposition}

\begin{proof}
Let $n \in \N$ and $x \in \R$, we introduce the (non-probability) measure $\P^*_x$ on $\hat{\calT}_n$, in which every possible choice of spine has mass $1$. More precisely, for any measurable function $f : \hat{\calT}_n \to \R_+$, we have $\int f d\P^*_x = \E_x\left[ \sum_{|w|=n} f(\T,V,w) \right]$. We compute by recurrence on $k \leq n$ the Radon-Nikod\'ym derivative of $\hat{\P}_x$ with respect to $\P^*_x$, to prove
\begin{equation}
  \label{eqn:rn1}
  \left.\frac{d\hat{\P}_x}{d\P^*_x}\right|_{\hat{\calF}_k} = \exp\left(\theta (V(w_k)-x) - \sum_{j=1}^k \kappa_j(\theta) \right).
\end{equation}
Observe that for $k=1$, \eqref{eqn:rn1} follows from the definition of $\hat{L}_1$ and $w(1)$. Writing $L_1$ a point process of law $\calL_1$ and $f$ a non-negative $\hat{\calF}_1$ measurable function,
\[
  \E\left[ f(\hat{L}_1, w(1)) \right] = \E\left[ \sum_{k=1}^{N_1} f(\hat{L}_1, k) \frac{e^{\theta \ell_1(k)}}{\sum_{j=1}^{N_k} e^{\theta \ell_1(j)}} \right]
  = \E\left[ \sum_{k=1}^{N_k} f(L_1,k)e^{\theta \ell_1(k) - \kappa_1(\theta)} \right].
\]
We now assume \eqref{eqn:rn1} true for some $k<n$, and we observe that
\begin{align*}
  \left.\frac{d\hat{\P}_x}{d\P^*_x}\right|_{\hat{\calF}_{k+1}}
  &= \left. \frac{d\hat{\P}_x}{d\P^*_x}\right|_{\hat{\calF}_k} \times \left(\sum_{u \in \Omega(w_k)} e^{\theta(V(u)-V(w_k)) - \kappa_{k+1}(\theta)}\right) \frac{e^{-V(w_{k+1})-V(w_k)}}{\sum_{u \in \Omega(w_k)} e^{\theta(V(u)-V(w_k)) - \kappa_{k+1}(\theta)}}\\
  &= \exp\left(\theta (V(w_k)-x) - \sum_{j=1}^k \kappa_j(\theta)\right)e^{\theta (V(w_{k+1})-V(w_k)) - \kappa_{k+1}(\theta)},
\end{align*}
which proves \eqref{eqn:rn1}.

As a consequence, for any $f : \calT \to \R_+$ measurable, we have
\begin{align*}
  \hat{\E}_x\left[ f( \T, V) \right]
  &= \int_{\hat{\calT}_n} e^{\theta (V(w)-x) - \sum_{j=1}^n \kappa_j(\theta)} f(\T,V) d\P_x^*(\T,V,w)\\
  &= \E_x\left[ f(\T,V) \sum_{|w|=n} e^{\theta (V(w)-x) - \sum_{j=1}^n \kappa_j(\theta)} \right]
  = e^{-\theta x} \E_x\left[ W_n f(\T,V) \right]
\end{align*}
therefore $\frac{\left.d\hat{\P}_x\right|_\calF}{d\P_x} = \frac{d\bar{\P}_x}{d\P_x} = e^{-\theta x} W_n$ which proves \eqref{eqn:spinaldecomposition}. Consequently, for any $f : \calT \to \R_+$ and $u \in \calU$ with $|u|=n$, we have
\begin{align*}
  \hat{\E}_x\left[ f(\T,V) \ind{w=u} \right]
  &= \int_{\hat{\calT}_n} e^{\theta (V(w)-x) - \sum_{j=1}^n \kappa_j(\theta)} f(\T,V) \ind{w=u} d \P^*_x(\T,V,w)\\
  &= \E_x\left[ f(\T,V) \sum_{|v|=n} e^{\theta (V(v)-x) - \sum_{j=1}^n \kappa_j(\theta)} \ind{v=u}\right]\\
  &= \E_x\left[ f(\T,V) e^{\theta (V(u) - x) - \sum_{j=1}^n \kappa_j(\theta)} \ind{u \in \T} \right]\\
  &= \hat{\E}_x \left[ \frac{e^{\theta (V(u)-x) - \sum_{j=1}^n \kappa_j(\theta)}}{e^{-\theta x}W_n} f(\T,V) \ind{u \in \T} \right]\\
  &= \hat{\E}_x \left[ \frac{e^{\theta V(u) - \sum_{j=1}^n \kappa_j(\theta)}}{W_n} f(\T,V) \ind{u \in \T} \right].
\end{align*}
\end{proof}

A direct consequence of this result, is the well-known many-to-one lemma. This equation, known at least from the early work of Peyrière \cite{Pey74} has been used in many forms over the last decades, and we introduce here its time-inhomogeneous version.

\begin{lemma}[Many-to-one]
\label{lem:manytoone}
We define an independent sequence of random variables $(X_k, k \leq n)$ that verifies
\[
  \forall k \leq n, \forall x \in \R, \P\left[ X_k \leq x \right] = \E\left[ \sum_{\ell \in L_k} \ind{\ell \leq x} ) e^{\theta \ell - \kappa_k(\theta)} \right].
\]
We write $S_k = S_0 + \sum_{j=1}^k X_j$ for $k \leq n$, where $\P_x(S_0=x)=1$. For all $x \in \R, k \leq n$ and measurable non-negative function $f$, we have
\begin{equation}
  \label{eqn:manytoone}
  \E_x\left[ \sum_{|u|=k} f(V(u_1),\ldots V(u_k)) \right] = e^{\theta x} \E_x\left[ e^{-\theta S_k + \sum_{j=1}^k \kappa_j(\theta)} f(S_1,\ldots S_k) \right].
\end{equation}
\end{lemma}

\begin{proof}
Let $f$ be a measurable non-negative function and $x \in \R$, we have, by Proposition \ref{prop:spinaldecomposition}
\begin{align*}
  \E_x \left[ \sum_{|u|=k} f(V(u_1),\ldots V(u_{k})) \right]
  &= \bar{\E}_{x}\left[ \frac{e^{\theta x}}{W_k}\sum_{|u|=k} f(V(u_1),\ldots V(u_k)) \right]\\
  &= \hat{\E}_{x}\left[ \frac{e^{\theta x}}{W_k}\sum_{|u|=k} f(V(u_1),\ldots V(u_k)) \right]\\
  &= \hat{\E}_{x}\left[ e^{-\theta (V(w_k)-x) + \sum_{j=1}^k \kappa_j(\theta)} f(V(w_1),\ldots V(w_k)) \right].
\end{align*}
We conclude noting that $(V(w_1), \ldots, V(w_n))$ under $\hat{\P}_x$ and $(S_1,\ldots S_n)$ under $\P_x$ have the same law.
\end{proof}

The many-to-one lemma and the spinal decomposition enable to compute additive moments of branching random walks, by using random walk estimates. These estimates are introduced in the next section, and extended to include time-inhomogeneous versions.

\section{Some random walk estimates}
\label{sec:randomwalk}

\subsection{Classical random walk estimates}
\label{subsec:classical}

We collect here well-known random walk estimates. We use these to compute the probability for a random walk with an interface to make an excursion above a given curve. The proofs, rather technical, are postponed to the Appendix \ref{app:estimates}.

\begin{figure}[ht]
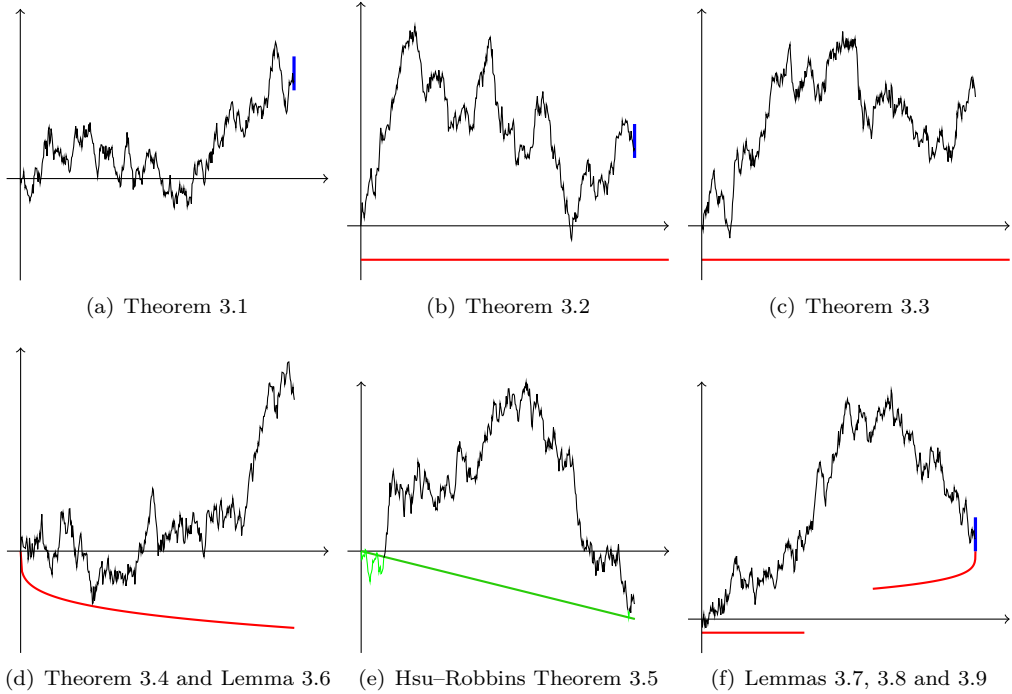

\label{fig:events}
\centering
\subfigure[Theorem \ref{thm:locallimit}]{

}%
\caption{Illustrations of the events controlled in Section \ref{subsec:classical}}
\end{figure}

We denote by $(T_n,n \geq 0)$ a one-dimensional centred random walk, with finite variance $\sigma^2$. The events we bound are illustrated in Figure \ref{fig:events}. We begin with a consequence of Stone's local limit theorem, which bounds the probability for a random walk to end up in an interval of finite size.
\begin{theorem}[Stone \cite{Sto65}]
\label{thm:locallimit}
There exists $C>0$ such that for all $a \geq 0$ and $h \geq 0$
\[ \limsup_{n \to +\infty} n^{1/2}\sup_{|y| \geq a n^{1/2}} \P(T_n \in [y,y+h]) \leq C (1+h)e^{-\frac{a^2}{2\sigma^2}}. \]
Moreover, there exists $H>0$ such that for all $a<b \in \R$
\[ \liminf_{n \to +\infty} n^{1/2}\inf_{y \in [an^{1/2},bn^{1/2}]} \P(T_n \in [y,y+H]) > 0. \]
\end{theorem}

Similar result is obtained by Caravenna and Chaumont, for a random walk conditioned to stay positive.
\begin{theorem}[Caravenna--Chaumont \cite{CaC08}]
\label{thm:locallimit_excursion}
Let $(r_n)$ be a positive sequence such that $r_n = O(n^{1/2})$. There exists $C>0$ such that for all $a \geq 0$ and $h \geq 0$,
\[ \limsup_{n \to +\infty} n^{1/2} \sup_{y \in [0,r_n]} \sup_{x \geq an^{1/2}} \P(T_n \in [x,x+h]|T_j \geq -y, j \leq n) \leq C (1 +h)ae^{-\frac{a^2}{2\sigma^2}}. \]
Moreover, there exists $H>0$ such that for all $a<b \in \R_+$,
\[ \liminf_{n \to +\infty} n^{1/2}\inf_{y \in [0,r_n]} \inf_{x \in [an^{1/2},bn^{1/2}]} \P(T_n \in [x,x+H]|T_j \geq -y, j \leq n) > 0. \]
\end{theorem}

Up to a transformation $T \mapsto T/(2H)$, which correspond to shrink the space by a factor $\frac{1}{2H}$, we assume in the rest of this article that random walks we consider are such that the lower bounds of Theorems \ref{thm:locallimit} and \ref{thm:locallimit_excursion} both hold with $H=1$.

The next result, often called in the literature the ``ballot theorem", give upper and lower bounds for the probability for a random walk to stay above zero. This result is stated in \cite{Koz76}, see also \cite{ABR08} for a review article on ballot theorems.
\begin{theorem}[Kozlov \cite{Koz76}]
\label{thm:ballot1}
There exists $C>0$ such that for all $n \geq 1$ and $y \geq 0$,
\[ \P(T_j \geq -y, j \leq n) \leq C (1 + y)n^{-1/2}. \]
Moreover, there exists $c>0$ such that for all $y \in [0,n^{1/2}]$
\[ \P(T_j \geq -y, j \leq n) \geq c(1 + y)n^{-1/2}.\]
\end{theorem}

A modification of this theorem, Theorem 3.2 of Pemantle and Peres in \cite{PeP95}, expresses the probability for a random walk to stay above some a boundary which moves ``strictly slower than $n^{1/2}$''.
\begin{theorem}[Pemantle--Peres \cite{PeP95}]
\label{thm:ballot}
Let $f: \N \to \N$ be an increasing positive function. The condition $\sum_{n \geq 0} \frac{f_n}{n^{3/2}}<+\infty$ is necessary and sufficient for the existence of an integer $n_f$ such that
\[ \sup_{n \in \N} n^{1/2} \P(T_j \geq -f_j, n_f \leq j \leq n) < +\infty. \]
\end{theorem}

To bound the probability for a random walk to stay above a linear boundary, we introduce the Hsu--Robbins theorem, bounding the expected number of times a random walk is above a linear boundary.
\begin{theorem}[Hsu--Robbins \cite{HsR47}]
\label{thm:hsurobbins}
For any $\epsilon>0$, we have $\sum_{n \geq 0} \P(T_n \leq -n\epsilon)<+\infty$.
\end{theorem}

We extend Theorems \ref{thm:ballot1} and \ref{thm:ballot}, to obtain a quantitative upper bound for the probability of a random walk to stay above a slowly moving boundary
\begin{lemma}
\label{lem:ballot}
Let $(f_n) \in \R^n$. If there exists $\alpha \in [0,1/2)$ and $A>0$ such that for any $n \in \N$, $|f_n| \leq A n^\alpha$ then there exists $C>0$ such that for all $y \geq 0$ and $n \geq 1$, we have
\[ \P(T_j \geq -y-f_j, j \leq n) \leq C(1 + y)n^{-1/2}. \]
\end{lemma}

In an earlier version of this article, we gave a wrong proof for this fact, as pointed out to us by Ming Fang. We present a corrected version, in the appendix. The proof uses results from \cite{DSV}.

The next lemma is also proved in Appendix \ref{app:estimates1}. The following upper bound of the probability for a random walk to make an excursion holds.
\begin{lemma}
\label{lem:excursion}
There exists $C>0$ such that for all $p,q \in \N$, $x,h \geq 0$  and $y \in \R$, we have
\begin{multline*}
  \P(T_{p+q} \in [y+h,y+h+1], T_j \geq -x \ind{j \leq p} + y \ind{p<j \leq p+q}, j \leq p+q)\\
  \leq C \frac{1+x}{p^{1/2}} \frac{1}{\max(p,q)^{1/2}} \frac{1+h}{q^{1/2}}.
\end{multline*}
\end{lemma}

We sum up Theorems \ref{thm:ballot1} and \ref{thm:ballot} and Lemmas \ref{lem:ballot} and \ref{lem:excursion} to obtain a general upper bound for the probability for a time-inhomogeneous random walk to make an excursion. Let $p,q,r \in \N$, we write $n=p+q+r$, $(X_k)_{k \in \N}$ and $(\tilde{X}_k)_{k \in \N}$ two independent families of i.i.d. random variables, with mean 0 and finite variance, and $(Y_n)_{n \geq 0}$ a family of independent random variables. We define the time-inhomogeneous random walk $(S_k, k \leq n)$ as follows:
\[ S_k = \sum_{j=1}^{\min\{k,p\}} X_j + \sum_{j=1}^{\min\{k-p,q\}}Y_j + \sum_{j=1}^{\min\{k-p-q,r\}}\tilde{X}_j. \]

Let $A \in \R$, and $x,y \in \R_+$, $h \in \R$, we denote by
\[\Gamma^{A,1}(x,y,h) = \{ s \in \R^n : \forall k \leq p, s_k \geq - x\}\]
the set of trajectories staying above $-x$ during the initial steps, and by
\[\Gamma^{A,3}(x,y,h) = \{ s \in \R^n : \forall k \in [n-r,n], s_k \geq y + A \log \tfrac{n}{n-k+1} \}. \]
\begin{lemma}
\label{lem:ballot_general}
For any $A \in \R$ and $F \subset \{1,3\}$,  there exists $C>0$ such that for all $p,q,r \in \N$, $x,y \in \R_+$ and $h \in \R$, we have
\begin{multline*}
  \P\left[ S_n + A \log n \in [y+h,y+h+1], (S_k, k \leq n) \in \bigcap_{f \in F} \Gamma^{A,f}(x,y,h)\right]\\
   \leq C \frac{1 + y \indset{F}(1)}{p^{\indset{F}(1)/2}} \frac{1}{\max(p,r)^{1/2}} \frac{1+h_+\indset{F}(3)}{r^{\indset{F}(3)/2}}.
\end{multline*}
\end{lemma}
This lemma is proved in Appendix \ref{app:estimates2}.

We finish this list of results with a lower bound of an event similar to the one studied in the previous lemma. We consider $\mu^{(1)}, \ldots \mu^{(P)}$ centred probability measures on $\R$ with finite variance. The process $S$ is defined as the sum of independent random variables such that for all $k \in [\alpha^{(n)}_{p-1}, \alpha^{(n)}_p)$, $S_{k+1}-S_k$ has law $\mu^{(p)}$. For $F \subset \{1,3\}$ and $x,y,\delta \in \R_+$, we write
\[
  \Upsilon^F(x,y,\delta) = \left\{s \in \R^n :
  \begin{array}{l}
  \forall k \leq \alpha_1^{(n)}, s_k \geq -x \ind{1 \in F} - \delta k \ind{1 \not \in F}\\
  \forall k \in (\alpha_1^{(n)}, \alpha_{P-1}^n], s_k \geq 0\\
  \forall k \in (\alpha_{P-1}^n, n], s_k \geq y \ind{3 \in F} -\delta(n-k)\ind{3 \not \in F}
  \end{array}
  \right\}.
\]
The next lemma bounds from below the probability for a random walk through a series of interfaces to be in $\Upsilon$. This lemma is proved in Appendix \ref{app:estimates3}.
\begin{lemma}
\label{lem:ballotlowerbound}
There exists $c>0$ such that for all $n\geq 1$ large enough, $F \subset \{1,3\}$, $x \in [0,n^{1/2}]$, $y \in [-n^{1/2},n^{1/2}]$ and $\delta>0$
\[
  \P(S_n \leq y + 1, S \in \Upsilon^F(x,y,\delta)) \geq c \frac{1 + x\indset{F}(1)}{n^{\indset{F}(1)/2}}\frac{1}{n^{1/2}}\frac{1}{n^{\indset{F}(3)/2}}
\]
\end{lemma}

\subsection{Extension to enriched random walks}

We extend here some of the results of the previous section to a random walk enriched with other random variables, which only depend on the last step of the random walk. We denote by $((X_n,\xi_n), n \geq 0)$ an i.i.d. sequence of random variables taking values in $\R^2$, such that $\E(X_1)=0$ and $\E(X_1^2)<+\infty$. We set $T_n = T_0 + X_1 + \cdots + X_n$, where $\P_x(T_0=x)=1$. The process $(T_n,\xi_n,n\geq 0)$ is an useful toy-model for the study of the spinal decomposition of the branching random walk, defined in Section 2.2. We begin with a lemma similar to Theorem \ref{thm:ballot1}.
\begin{lemma}
\label{lem:ballot_spine}
We suppose that $\E(X_1)=0$, $\E(X_1^2)<+\infty$ and $\E((\xi_1)_+^2)<+\infty$. There exists $C>0$ that does not depend on the law of $\xi_1$ such that for any $n \in \N$ and $x \geq 0$, we have
\[ \P_x\left[ T_j \geq 0, j \leq n, \exists k \leq n : T_k \leq \xi_k \right] \leq C \frac{1+x}{n^{1/2}}\left[\P(\xi_1 \geq 0) + \E((\xi_1)_+^2)\right]. \]
\end{lemma}

\begin{proof}
Let $n \in \N$ and $x \geq 0$. We observe that
\[ \P_x\left[ T_j \geq 0, j \leq n, \exists k \leq n : T_k \leq \xi_k \right] \leq \sum_{k=1}^n \underbrace{\P_x\left[T_k \leq \xi_k, T_j \geq 0, j \leq n\right]}_{\pi_k}. \]
Applying the Markov property at time $k$, we obtain
\[ \pi_k \leq \E_x\left[ \ind{T_k \leq \xi_k} \ind{T_j \geq 0, j \leq k} \P_{T_k}\left( T_j \geq 0, j \leq n-k \right)  \right]. \]
By use of Theorem \ref{thm:ballot1}, for all $z \in \R$, we have
\[ \P_z \left[ T_j \geq 0, j \leq n-k \right] \leq C (1 + z) (n-k+1)^{-1/2}\ind{z\geq 0}. \]
Thus, writing $(X,\xi)$ for a copy of $(X_1, \xi_1)$ independent of $(T_n, \xi_n, n \geq 0)$, we have
\begin{align*}
  \pi_k &\leq C (n-k+1)^{-1/2} \E_x\left[ \ind{\xi_k \geq 0} (1 + \xi_k) \ind{T_k \leq \xi_k} \ind{T_j \geq 0, j \leq k} \right] \\
  &\leq C (n-k+1)^{-1/2} \E_x \left[ \ind{\xi \geq 0} (1 + \xi_+) \ind{T_{k-1} \leq \xi_+  + X_-} \ind{T_j \geq 0, j \leq (k-1)} \right]
\end{align*}
We bound this quantity by conditioning on the value $\zeta = \xi_+ + X_- \geq 0$, we obtain
\[
\P_x(T_k \leq \zeta, T_j \geq 0, j \leq k) \leq 
  \begin{cases}
    C \frac{(1+x)(1 + \zeta^2)}{(k+1)^{3/2}} & \text{if } \zeta^2 \leq k, \text{ by Lemma \ref{lem:excursion}}\\
    C \frac{1+x}{(k+1)^{1/2}} & \mathrm{otherwise, by Theorem \ref{thm:ballot1}} 
  \end{cases}
\]
Summing all these estimates, we obtain
\begin{align*}
  \sum_{k=0}^{n-1} \frac{ \P_x \left[ T_k \leq \zeta, T_j \geq 0, j \leq k \right]}{(n-k+1)^{1/2}}
  \leq & C (1 + x) \sum_{k=0}^{\min(\zeta^2,n-1)} \frac{1}{(n-k+1)^{1/2}(k+1)^{1/2}}\\
  &\qquad \qquad + C (1 + x)(1+ \zeta^2)\sum_{k=\zeta^2}^{n-1} \frac{1}{(k+1)^{3/2}(n-k+1)^{1/2}} \\
  \leq & C(1 + x) (1 + \zeta) n^{-1/2}.
\end{align*}

As a consequence,
\begin{multline*}
  \sum_{k=1}^n \pi_k \leq C \frac{1 + x}{n^{1/2}} \E\left[ \ind{\xi \geq 0} (1 + X_- + \xi_+)(1+\xi_+) \right]\\
  \leq C \frac{1+x}{n^{1/2}}\left[ 1 + \E\left(X_-^2\right) \right] \left[ \P(\xi \geq 0) + \E\left(\xi_+^2\right)\right]
\end{multline*}
by Cauchy-Schwarz estimate, which ends the proof.
\end{proof}

We continue by reprising Lemma \ref{lem:excursion}. 
\begin{lemma}
\label{lem:excursion_spine}
We assume that $\E(X_1)=0$, $\E(X_1^2)<+\infty$ and $\E((\xi_1)_+^2)<+\infty$. For any $t \in (0,1)$, there exists $C>0$ that does not depend of the law of $\xi_1$, such that for all $n \in \N$, $x, h\geq 0$ and $y \in \R$, we have
\begin{multline*}
  \P_x\left[ T_n-y-h \in [0,1], T_j \geq y \ind{j>tn}, j \leq n, \exists k \leq n : T_k \leq \xi_k + y\ind{k>tn} \right]\\
  \leq C \frac{(1+x)(1+h)}{n^{3/2}} \left[\P(\xi_1 \geq 0) + \E((\xi_1)_+^2)\right].
\end{multline*}
\end{lemma}

\begin{proof}
Let $n \in \N$, $x,h \geq 0$ and $y\in \R$. We denote by $p=\floor{tn}$ and by 
\[\tau = \inf\{ k \geq 0 : T_k \leq \xi_k + y \ind{k>p} \}.\]
We observe that
\begin{multline}
  \label{eqn:itsall}
  \P_x\left[ T_n-y-h \in [0,1], T_j \geq y \ind{j>tn}, j \leq n, \tau \leq n \right]\\
  \leq \P_x\left(T_n-y-h \in [0,1], T_j \geq y \ind{j>p}, \tau \leq p\right)\\
  + \P_x\left(T_n-y-h \in [0,1], T_j \geq y \ind{j>p}, p < \tau \leq n\right).
\end{multline}

We first take interest in the event $\{\tau \leq p\}$. Applying the Markov property at time $p$, we obtain
\begin{equation}
  \label{eqn:partone}
  \P_x\left[T_n-y-h \in [0,1], T_j \geq y \ind{j>p}, \tau \leq p \right] = \E_x\left[\ind{T_j \geq 0, j \leq p} \ind{\tau \leq p} \phi(T_p)\right],
\end{equation}
writing $\phi(z) = \P_z \left[ T_{n-p}-y-h \in [0,1], T_j \geq y, j \leq n-p\right]$, for $z \in \R$. Applying Lemma \ref{lem:ballot_general}, we have $\sup_{z \in \R} \phi(z) \leq C (1 +h)n^{-1}$. Therefore
\begin{align*}
  \P_x\left(T_n-y-h \in [0,1], T_j \geq y \ind{j>p}, \tau \leq p\right)
  &\leq C \frac{1 +h}{n} \P_x\left[ T_j \geq 0, j \leq p, \exists k \leq p : T_k \leq \xi_k \right]\\
  &\leq C \frac{(1 +x)(1+h)}{n^{3/2}} \left[ \P(\xi_1 >0) + \E\left((\xi_1)_+^2\right) \right]
\end{align*}
by use of Lemma \ref{lem:ballot_spine}.

We now take care of $\{ \tau > p\}$. We have
\begin{align*}
  &\P_x\left[ T_n-y-h \in [0,1], T_j \geq y \ind{j>p}, j \leq n, p \leq \tau \leq n \right]\\
  &\qquad\qquad\qquad\qquad \leq \P_x\left[
  \begin{array}{l}
     T_n-y-h \in [0, 1], T_n - T_{n-j} \leq y+h+1 - y \ind{n-j <p}\\ \exists k \leq n-p : T_{n-k} \leq \xi_{n-k} + y
  \end{array}
  \right]\\
  &\qquad\qquad\qquad\qquad \leq \P_x \left[
  \begin{array}{l}
     T_n - T_0-y-h+x \in [0,1], T_n - T_{n-j} \leq h+1+y \ind{j\geq n-p}\\ \exists k \leq n-p : T_n - T_{n-k} \geq y+h - (\xi_{n-k}+y)
  \end{array}
  \right].
\end{align*}
We denote by $\hat{T}_j = T_{n}-T_{n-j}$ and $\hat{\xi}_j=\xi_{n-j}$, we have
\begin{align*}
  &\P_x\left[ T_n-y-h \in [0,1], T_j \geq y \ind{j>p}, j \leq n, p \leq \tau \leq n \right]\\
  \leq &\P_x\left[ \hat{T}_n-y-h+x \in [0, 1], \hat{T}_j \leq h+1-y\ind{j \geq n-p}, \exists k \leq n-p : \hat{T}_k \geq y+h -(\hat{\xi}_{k} + y)\right].
\end{align*}
We observe that $(\hat{T}_j,\hat{\xi}_j, j \leq n)$ has the same law as $(T_j,\xi_j,j \leq n)$ under $\P_0$, as a consequence
\begin{align*}
  &\P_x\left[ T_n-y-h \in [0,1], T_j \geq y \ind{j>p}, j \leq n, p \leq \tau \leq n \right]\\
  \leq & \P_0 \left[ T_n-h-y+x \in [0,1], T_j \leq h+1 - y\ind{j \geq n-p}, \exists k \leq n-p : T_k \geq h- \xi_k \right]\\
  \leq & \P_{-h-1} \left[ T_n-y+x \in [-1,0], T_j \leq -y \ind{j \leq n-p} \exists k \leq n-p : T_k \geq -\xi_k - 1\right].
\end{align*}
This quantity is bounded in \eqref{eqn:partone}, replacing $(T,\xi)$ by $(-T,\xi)$, and exchanging the roles of $x$ and $h$, thus similar computations lead to
\begin{multline}
  \label{eqn:parttwo}
  \P_x\left[ T_n-y-h \in [0,1], T_j \geq y \ind{j>p}, j \leq n,p< \tau \leq n\right]\\ \leq C \frac{(1+x)(1+h)}{(n+1)^{3/2}} \left[\P(\xi \geq 0) + \E(\xi_+^2)\right]
\end{multline}
which ends the proof.
\end{proof}

We end with an analogue of the Hsu--Robbins theorem.
\begin{lemma}
\label{lem:hsurobbins_spine}
We suppose that $\E(X_1)=0$, $\E(X_1^2)<+\infty$ and $\E((\xi_1)_+)<+\infty$. Let $\epsilon>0$, there exists $C>0$ that does not depend on the law of $\xi_1$ such that for all $x,z \geq 0$ and $n \in \N$
\begin{equation*}
  \P_x\left[ T_j \geq -\epsilon j, j \leq n, \exists k \leq n : T_k \leq -\epsilon k + \xi_k \right] \leq C \left[\frac{\E\left[ (\xi+z)_+\right]}{\epsilon}\right] + \E_z\left[ \sum_{n \geq 0} \ind{T_n \leq -n \epsilon/2} \right].
\end{equation*}
\end{lemma}

\begin{proof}
By union bound, we have
\[
  \P_x\left[ T_j \geq - \epsilon j, j \leq n, \exists k \leq n : T_k \leq -\epsilon k + \xi_k \right] \leq \sum_{k=1}^n \P_x(T_k \leq -\epsilon k + \xi_k).
\]
Moreover $\P_x(T_k \leq - \epsilon k + \xi_k) \leq \P_x(T_k \leq z-\epsilon k / 2) + \P(\xi_k \geq \epsilon k / 2 + z)$, thus
\begin{align*}
  \P_x\left[\begin{array}{l} T_j \geq - \epsilon j, j \leq n\\ \exists k \leq n : T_k \leq -\epsilon k + \xi_k\end{array} \right]
  &\leq \sum_{k=1}^n \P_x\left[ T_k \leq -\epsilon j/2 - z\right] + \sum_{k=1}^n \P(\xi_k \geq \epsilon k/2 - z)\\
  &\leq \E_{x+z}\left[ \sum_{k =1}^{+\infty} \ind{T_j \leq -\epsilon j/2} \right] + 2 \frac{\E((\xi+z)_+)}{\epsilon}.
\end{align*}
\end{proof}

\begin{remark}
By dominated convergence theorem and Theorem \ref{thm:hsurobbins},
\[\lim_{z \to +\infty} \E_z[ \sum_{n \geq 0} \ind{T_n \leq -n \epsilon/2}] = 0,\]
thus to obtain a good bound in Lemma \ref{lem:hsurobbins_spine}, it is useful to choose $z$ very large.
\end{remark}

\section{Bounds on the tail of the maximal displacement}
\label{sec:tailbehaviour}

Let $(\T,V)$ be a BRWis of length $n$. We recall that $(\calL_p, p \leq P)$ is a family of point processes, and $0=\alpha_0 < \alpha_1 < \cdots < \alpha_P=1$ a sequence of real numbers. Up to replacing these sequences with $(\calL_1,\calL_1,\calL_1,\calL_2, \ldots \calL_P)$ and $0=\alpha_0<\alpha_1/3<2\alpha_1/3 < \alpha_1< \alpha_2<\ldots < \alpha_P=1$, we assume that $P \geq 3$. For $p \leq P$ and $\theta>0$, we write $\kappa_p(\theta) = \log \E\left[ \sum_{\ell \in L_p} e^{\theta \ell} \right]$ the log-Laplace transform of $\calL_p$. We write $M_n$ the maximal displacement at time $n$ of the BRWis. The main goal of this section is to prove the following estimate. 
\begin{theorem}
\label{thm:tailestimate}
Under the assumptions \eqref{eqn:optimization_problem}, \eqref{eqn:differentiable} and \eqref{eqn:variance}, there exists $C>0$ such that for all $n \in \N$ and $y \geq 0$,
\[
  \P(M_n \geq nv_\mathrm{is} - \lambda \log n + y) \leq C (1 + y\indset{B}(1)) e^{-\theta_1 y}.
\]
where $v_\mathrm{is}$ and $\lambda$ are defined respectively by \eqref{eqn:speedisdef} and \eqref{eqn:logisdef}. Moreover, under the additional assumption \eqref{eqn:integrability}, there exists $c>0$ such that for any $n \in \N$ large enough and $y \in [0,n^{1/2}]$, we have
\[
   \P(M_n \geq nv_\mathrm{is} - \lambda \log n + y) \geq c (1 + y\indset{B}(1))e^{-\theta_1 y}.
\]
\end{theorem}
To prove this result, we use the decomposition of the BRWis obtained thanks to Proposition \ref{prop:optimization_problem}. According to this result, if $\bba$ is the solution of \eqref{eqn:optimization_problem}, and $\theta_p=(\kappa^*_p)'(a_p)$, the sequence $\bbtheta$ is non-decreasing, and takes a finite number $T$ of values. We prove Theorem \ref{thm:tailestimate} by induction on $T$. In the next section, we prove Theorem \ref{thm:tailestimate} for a BRWis such that $T=1$. In Section \ref{subsec:extension}, we prove the induction hypothesis. Section \ref{subsec:proof} derives Theorem \ref{thm:main} from Theorem \ref{thm:tailestimate}.

\subsection{The case of a mono-parameter branching random walk}
\label{subsec:monoparameter}

We consider in a first time a BRWis $(\T,V)$ satisfying additional assumptions that guarantee the sequence $\bbtheta$ to be constant. We write,
\[ \forall \phi \in \R_+, \forall p \leq P, \quad E_p(\phi) = \sum_{q=1}^p (\alpha_{q}-\alpha_{q-1})(\phi \kappa'_q(\phi)-\kappa_q(\phi)).\]
We assume there exists $\theta > 0$ such that
\begin{equation}
  \label{eqn:monooptimalpath}
  \forall p \leq P,  E_p(\theta) \leq 0 \quad \mathrm{and} \quad E_P(\theta)=0.
\end{equation}
We write $a_p = \kappa'_p(\theta)$ and $B=\{p \leq P : E_p(\theta) = E_{p-1}(\theta)=0\}$. By \eqref{eqn:legendreestimate} and \eqref{eqn:monooptimalpath}, $\bba \in \calR$, and by Proposition \ref{prop:optimization_problem}, $\bba$ is the solution of \eqref{eqn:optimization_problem}. With these notations, we have
\begin{equation}
  \label{eqn:speedcomputed}
  v_\mathrm{is} = \sum_{p=1}^P (\alpha_p-\alpha_{p-1}) a_p \quad \mathrm{and} \quad \lambda = \frac{1}{2\theta}  \left( 1 + \indset{B}(1) + \indset{B}(P)\right).
\end{equation}

\begin{theorem}
  \label{thm:light}
Under assumptions \eqref{eqn:variance} and \eqref{eqn:monooptimalpath}, there exists $C>0$ such that for any $n \in \N$ and $y \geq 0$, we have
\[
  \P(M_n \geq nv_\mathrm{is} - \lambda \log n + y) \leq C (1 + y\indset{B}(1))e^{-\theta y}.
\]  
Moreover, under the additional assumption \eqref{eqn:integrability}, there exists $c>0$ such that for any $n \in \N$ large enough and $y \in [0,\sqrt{n}]$,
  \[ \P(M_n \geq nv_\mathrm{is} - \lambda \log n + y) \geq c (1 + y\indset{B}(1))e^{-\theta y}. \]
\end{theorem}

We write $m_n = n v_\mathrm{is} - \lambda \log n$ the expected position of $M_n$. To obtain the upper bound, we prove in a first time that with high probability, if the optimal path stays close to the boundary of the branching random walk during the first or the last time interval, then there is no individual above this boundary at any time. In a second time, we bound from above and from below the number of individuals who stayed below the boundary, and end at time $n$ close to $m_n$.

We introduce
\[K^{(n)}_k = \sum_{p=1}^P  \sum_{j=1}^k \kappa_p(\theta) \ind{j \in [\alpha^{(n)}_{p-1},\alpha^{(n)}_p)} \quad \mathrm{and}\quad \bar{a}^{(n)}_k = \sum_{p=1}^P  \sum_{j=1}^k a_p \ind{j \in [\alpha^{(n)}_{p-1},\alpha^{(n)}_p)}.\]
Using Equation \eqref{eqn:legendreestimate}, we observe that
\begin{equation}
  \label{eqn:monolegendre} \theta \bar{a}^{(n)}_k - K^{(n)}_k = \sum_{p=1}^P \kappa^*_p(a_p) \sum_{j=1}^k \ind{j \in [\alpha^{(n)}_{p-1},\alpha^{(n)}_p)},
\end{equation}
thus, writing $\phi_t = \int_0^t \sum_{p=1}^P \kappa^*_p(a_p) \ind{s \in [\alpha_{p-1}, \alpha_p)}ds$, we have
\begin{equation}
  \label{eqn:monoestimate}
  \sup_{n \geq 0} \sup_{k \leq n} \left| \theta \bar{a}^{(n)}_k - K^{(n)}_k - n\phi_\frac{k}{n} \right| < +\infty.
\end{equation}

\subsubsection{A frontier for the branching random walk}
\label{subsec:upperbound}

We prove in a first time that if $1 \in B$, then with high probability, there is no individual to the right of $a_1 k$ at any time $k \leq \alpha^{(n)}_1$.
\begin{lemma}
\label{lem:upperfrontier_starter}
Under assumption \eqref{eqn:monooptimalpath}, if $1 \in B$, then for all $y \geq 0$ and $n \in \N$,
\[ \P( \exists u \in \T, |u| \leq \alpha^{(n)}_1 : V(u) \geq a_1 |u| + y) \leq e^{-\theta y}. \]
\end{lemma}

\begin{proof}
Let $y \geq 0$ and $n \geq 1$. For $k \leq \alpha^{(n)}_1$, we write
\[Z^{(n)}_k = \sum_{|u|=k} \ind{V(u) \geq a_1 k + y} \ind{V(u_j) \leq a_1 j +y, j \leq k}\]
the number of individuals for the first time at time $k$ above the curve $a_1 \cdot + y$. By use of \eqref{eqn:manytoone}, we have
\[
  \E(Z^{(n)}_k) = \E\left[ e^{-\theta S_k + k \kappa_1(\theta)} \ind{S_k \geq k a_1 + y} \ind{S_j \leq j a_1 + y, j < k} \right], \]
where $S$ is a random walk with mean $\E\left[ \sum_{\ell \in L_1} \ell e^{\theta \ell - \kappa_1(\theta)} \right] = \kappa'_1(\theta)=a_1$
and finite variance. Moreover, as $1 \in B$, we have $E_1=\theta a_1 - \kappa_1(\theta)=0$ and
\[
  \E(Z^{(n)}_k) \leq e^{-\theta y} \P(S_k \geq ka_1 + y, S_j \leq j a_1 +y,j<k).
\]
As a consequence, by Markov inequality, we have
\begin{multline*}
  \P( \exists u \in \T, |u| \leq \alpha^{(n)}_1 : V(u) \geq a_1 |u| + y) \leq \sum_{k=1}^{\alpha^{(n)}_1} \E(Z^{(n)}_k)\\
  \leq e^{-\theta y} \sum_{k=1}^n \P(S_k \geq ka_1 + y, S_j \leq j a_1 +y, j < k) \leq e^{-\theta y} \P(\exists k \leq n : S_k \geq k a_1 + y).
\end{multline*}
\end{proof}

We compute, if $P \in B$, the probability that there exists at some time $k \geq \alpha^{(n)}_{P-1}$ an individual above some well-chosen curve. To do so, we denote by $r^{(n)}_k = a_P(k-n) + \frac{3}{2\theta} \log (n-k+1)$. We add a piece of notation to describe the frontier of the branching random walk. We write
\[
  F^{(n)} = \bigcup_{p \in B \cap \{1,P\}} \left[\alpha^{(n)}_{p-1}, \alpha^{(n)}_p\right], \quad F^{(n)}_k = F^{(n)} \cap [0,k]\]
and $f^{(n)}_j = a_1 j \ind{j \leq \alpha^{(n)}_1} + ( m_n + r^{(n)}_j ) \ind{j \geq \alpha^{(n)}_{P-1}}$ for any $j \in F^{(n)}$. The following estimate holds.
\begin{lemma}
\label{lem:upperfrontier_end}
Under assumptions \eqref{eqn:variance} and \eqref{eqn:monooptimalpath}, if $P \in B$, there exists $C>0$ such that for all $y \geq 0$ and $n \in \N$,
\[
 \P\left[ \exists |u| > \alpha^{(n)}_{P-1} : V(u) \geq m_n + r^{(n)}_k +y \right]  \leq C(1 +y\indset{B}(1)) e^{-\theta y}.
\]
\end{lemma}

\begin{proof}
We assume in a first time that $1 \not \in B$. We have $\lambda = \frac{1}{\theta}$ and
\begin{align}
  &\P\left[ \exists |u| > \alpha^{(n)}_{P-1} : V(u) \geq  m_n + r^{(n)}_k +y \right] \nonumber\\
  \leq &\E\left[ \sum_{|u| \geq \alpha^{(n)}_{P-1}} \ind{V(u) \geq m_n + r^{(n)}_{|u|} + y} \ind{V(u_j) \leq m_n + r^{(n)}_j + y, \alpha^{(n)}_{P-1} \leq j < k} \right] \nonumber\\
  \leq &\sum_{k=\alpha^{(n)}_{P-1}}^n \E\left[ e^{-\theta S_k + K^{(n)}_k} \ind{S_k \geq m_n + r^{(n)}_k + y, S_j \leq m_n + r^{(n)}_j + y, \alpha^{(n)}_{P-1} \leq j < k} \right]\nonumber\\
  \leq &C\sum_{k=\alpha^{(n)}_{P-1}}^n \frac{n^{\theta \lambda} e^{-\theta y}}{(n-k+1)^{3/2}} \P\left( S_k \geq m_n + r^{(n)}_k + y, S_j \leq m_n + r^{(n)}_j + y, \alpha^{(n)}_{P-1} \leq j < k \right), \label{eqn:coucou}
\end{align}
by \eqref{eqn:manytoone} and \eqref{eqn:monoestimate}. By conditioning with respect to $S_k - S_{k-1}$, we have
\[  \P\left( S_k \geq m_n + r^{(n)}_k + y, S_j \leq m_n + r^{(n)}_j + y, \alpha^{(n)}_{P-1} \leq j < k \right)
  = \E\left[ \phi_k(S_k - S_{k-1} - a_1) \right],
\]
writing for $x \in \R$,
\begin{align*}
  \phi_k(x) &= \P(S_{k-1} \geq m_n + r^{(n)}_k + y-x, S_j \leq m_n + r^{(n)}_j + y, \alpha^{(n)}_{P-1} \leq j \leq k-1)\\
  &= \sum_{h=0}^{+\infty} \P\left(
  \begin{array}{l}
    S_j \leq m_n + r^{(n)}_j + y, \alpha^{(n)}_{P-1} \leq j \leq k-1 \\ S_{k-1} - m_n - r^{(n)}_k - y - h \in [h, h + 1)
  \end{array}
  \right)\\
  &\leq \sum_{h=0}^\floor{x} C \frac{1+h}{n^{1/2}(k-\alpha^{(n)}_{P-1})^{1/2}} \leq C \frac{(1 + x_+)^2}{n^{1/2}(k-\alpha^{(n)}_{P-1})^{1/2}},
\end{align*}
by Lemma \ref{lem:ballot_general}. We have
\[  \P\left( S_k \geq m_n + r^{(n)}_k + y, S_j \leq m_n + r^{(n)}_j + y, \alpha^{(n)}_{P-1} \leq j < k \right) \leq \frac{C}{n^{1/2}(k - \alpha^{(n)}_{P-1})^{1/2}}, \]
as a consequence \eqref{eqn:coucou} becomes
\begin{multline*}
  \P\left[ \exists |u| > \alpha^{(n)}_{P-1} : V(u) \geq  m_n + r^{(n)}_k +y \right]\\
  \leq \sum_{k=\alpha^{(n)}_{P-1}}^n C e^{-\theta y} \frac{n^{1/2}}{(k-\alpha^{(n)}_{P-1}+1)^{1/2}(n-k+1)^{3/2}}
  \leq C e^{-\theta y}.
\end{multline*}

In a second time, if $1 \in B$, then $\lambda = \frac{3}{2\theta}$. We have
\begin{multline*}
  \P\left[ \exists |u| > \alpha^{(n)}_{P-1} : V(u) \geq  m_n + r^{(n)}_k +y \right] \\
  \leq \P\left[ \exists |u| \leq \alpha^{(n)}_1 : V(u) \geq a_1 |u| + y \right] \qquad \qquad \qquad \qquad \qquad \qquad \qquad \qquad \\
  + \P\left[ \exists |u| \geq \alpha^{(n)}_{P-1} : V(u) \geq f^{(n)}_k + y, V(u_j) \leq f^{(n)}_j + y, j \in F^{(n)}_{k-1} \right].
\end{multline*}
Lemma \ref{lem:upperfrontier_starter} bounds the first part of this inequality. By \eqref{eqn:manytoone}, for $k \geq \alpha^{(n)}_{P-1}$ we have
\begin{align*}
  & \E\left[ \sum_{|u|=k} \ind{V(u) \geq f^{(n)}_k  + y} \ind{V(u_j) \leq f^{(n)}_j +y, j \in F^{(n)}_{k-1}} \right]\\
  &\qquad \qquad \leq \E\left[ e^{-\theta S_k + K^{(n)}_k} \ind{S_k \geq f^{(n)}_k+ y} \ind{S_j \leq f^{(n)}_j +y, j \in F^{(n)}_{k-1}} \right]\\
  &\qquad \qquad \leq C \frac{n^{\theta \lambda}}{(n-k+1)^{3/2}}e^{- \theta y} \P(S_k \geq f^{(n)}_k +y, S_j \leq f^{(n)}_j +y, j \in F^{(n)}_{k-1})\\
  &\qquad \qquad \leq C (1 + y) e^{-\theta y} \frac{n^{3/2}}{(k - \alpha^{(n)}_{P-1}+1)^{3/2} (n-k+1)^{3/2}},
\end{align*}
using again Lemma \ref{lem:ballot_general}, and conditioning with respect to the last step of the random walk. By Markov inequality, we have
\begin{multline*}
  \P\left[ \exists |u| \geq \alpha^{(n)}_{P-1} : V(u) \geq f^{(n)}_k + y, V(u_j) \leq f^{(n)}_j + y, j \in F^{(n)}_{k-1} \right]\\
  \leq  C (1+y)e^{-\theta y} \sum_{k=\alpha^{(n)}_{P-1}}^n \frac{n^{3/2}}{(k - \alpha^{(n)}_{P-1}+1)^{3/2} (n-k+1)^{3/2}}
  \leq C (1+y)e^{-\theta y},
\end{multline*}
ending the proof.
\end{proof}

These two lemmas imply that with high probability, there is no individual above $f^{(n)}+y$ at any time in $F^{(n)}$. To complete the proof of the upper bound for the tail distribution of $M_n$, we compute the number of individuals who, travelling below that boundary, are at time $n$ in a neighbourhood of $m_n$. We write
\[ X^{(n)}(y,h) = \sum_{|u| = n} \ind{V(u) - m_n - y \in [-h,-h+1]}\ind{V(u_j) \leq f^{(n)}_j+y, j \in F^{(n)}}. \]
\begin{lemma}
\label{lem:upperestimate}
Under assumptions \eqref{eqn:variance} and \eqref{eqn:monooptimalpath}, there exists $C>0$ such that for all $n \geq 1$, $y \in \R_+$ and $h \in \R$, we have
\[ \E(X^{(n)}(y,h)) \leq C (1 + y\indset{B}(1))(1 + h_+ \indset{B}(P)) e^{-\theta(y-h)}. \]
\end{lemma}

\begin{proof}
Note that if $P \in B$ and $h<-1$, then $X^{(n)}(y,h)=0$. Otherwise, using Equation \eqref{eqn:manytoone}, we have
\begin{align*}
  \E(X^{(n)}(y,h))
  &=\E\left[ e^{-\theta S_n + K^{(n)}_n} \ind{S_n - m_n - y \in [-h,-h+1]} \ind{S_j \leq f^{(n)}_j+y, j \in F^{(n)}} \right]\\
  &\leq C n^{\theta \lambda}e^{-\theta (y-h)} \P\left(S_n - m_n - y \in [-h,-h+1], S_j \leq f^{(n)}_j+y, j \in F^{(n)}\right)
\end{align*}
by Equation \ref{eqn:monoestimate}. Applying Lemma \ref{lem:ballot_general}, we obtain
\[
   \P\left(S_n - f^{(n)}_n - y \in [-h,-h+1], S_j \leq f^{(n)}_j+y, j \in F^{(n)}\right)
   \leq C\frac{(1 + y\indset{B}(1))(1 + h_+ \indset{B}(P))}{(n+1)^{(1+\indset{B}(1)+\indset{B}(P))/2}}.
\]
\end{proof}

These lemmas can be used to obtain a tight upper bound for $\P(M_n \geq nv_\mathrm{is}-\lambda \log n + y)$.
\begin{corollary}
Under assumptions \eqref{eqn:variance} and \eqref{eqn:monooptimalpath}, there exists $C>0$ such that for all $y \geq 0$ and $n \in \N$, we have
\[ \P(M_n \geq nv_\mathrm{is} - \lambda \log n +y) \leq C(1 + y\indset{B}(1))e^{-\theta y}. \]
\end{corollary}

\begin{proof}
Let $y \geq 0$ and $n \in \N$, we have
\[
  \P(M_n \geq nv_\mathrm{is} - \lambda \log n +y)
  \leq \P\left(\exists |u| \in F^{(n)} : V(u) \geq f^{(n)}_{|u|}+y\right)
  + \sum_{h=0}^{+\infty} \E\left(X^{(n)}(y,-h)\right).
\]
Using Lemmas \ref{lem:upperfrontier_starter} and \ref{lem:upperfrontier_end}, we have $\P(\exists |u| \in F^{(n)} : V(u) \geq f^{(n)}_{|u|}+y) \leq C(1 +y\indset{B}(1))e^{-\theta y}$. Applying Lemma \ref{lem:upperestimate}, we obtain
\[ \sum_{h=0}^{+\infty} \E(X^{(n)}(y,-h)) \leq C (1 +y\indset{B}(1))e^{-\theta y} \sum_{h=0}^{+\infty} e^{-\theta h} \leq C (1 +y\indset{B}(1))e^{-\theta y}. \]
\end{proof}

\subsubsection{Lower bound through a second order computation}
\label{subsec:lowerbound}

To bound from below $\P(M_n \geq m_n+y)$, we bound from below the probability there exists an individual alive at time $n$, which stayed an any time $k \leq n$ below some curve $g^{(n)}$ defined below and is at time $n$ above $m_n$. We write $B^{(n)} = \cup_{p \in B} (\alpha^{(n)}_{p-1}, \alpha^{(n)}_p]$ the set of times such that the optimal path is close to the frontier of the BRWis. We choose $\delta>0$ small enough such that $3 \theta \delta < \min_{p \in B^c} -E_p(\theta)$. For all $n \geq 1$, $p \leq P$ and $k \in (\alpha^{(n)}_{p-1},\alpha^{(n)}_p]$ we define
\[ g^{(n)}_k = 1+
\begin{cases}
  \bar{a}^{(n)}_k - \ind{p=P} \lambda \log n & \mathrm{if} \quad E_p(\theta)=E_{p-1}(\theta)=0\\
  \bar{a}^{(n)}_k + (k-\alpha^{(n)}_{p-1}) \delta & \mathrm{if} \quad E_{p-1}(\theta)=0,E_p(\theta)<0\\
  \bar{a}^{(n)}_k + (\alpha^{(n)}_p-k) \delta & \mathrm{if} \quad E_p(\theta)=0, E_{p-1}(\theta)<0\\
  \bar{a}^{(n)}_k + \delta n &\mathrm{otherwise.}
\end{cases} \]
With this definition, using \eqref{eqn:monoestimate}, we have,
\begin{equation}
  \label{eqn:polyestimate}
  \theta g^{(n)}_{k}-K^{(n)}_{k} \leq C +
  \begin{cases}
    -\ind{p=P} \theta \lambda \log n  & \quad E_p(\theta)=E_{p-1}(\theta)=0\\
    -\delta (k-\alpha^{(n)}_{p-1}) & \mathrm{if} \quad E_{p-1}(\theta)=0, E_p(\theta)<0\\
    -\delta (\alpha^{(n)}_{p+1}-k)  - \ind{p=P} \theta \lambda \log n &\mathrm{if} \quad E_{p-1}(\theta)<0, E_p(\theta)=0\\
    -\delta n & \mathrm{if} \quad E_{p-1}(\theta)>0, E_p(\theta)>0.
  \end{cases}
\end{equation}

We prove in the rest of the section that the set
\[ \calA_n(y) = \left\{ u \in \T : V(u) \geq m_n + y, V(u_j) \leq g^{(n)}_j + y, j \leq n  \right\} \]
is non-empty. To do so, we restrict this set to individuals with a constraint on their reproduction. For $u \in \T$, we denote by
\[ \xi(u) = \sum_{u' \in \Omega(u)} \left(1 + (V(u')-V(u))_+\ind{|u| \in B^{(n)}+1}\right) e^{\theta (V(u')-V(u))} \]
a quantity closely related to the spread of the offspring of $u$. We write, for $z>0$ and $p \leq P$
\[ \calB_n(z) = \left\{ u \in \T : |u| = n, \xi(u_j) \leq z e^{-\frac{\theta}{2}\left[ V(u_j)-g^{(n)}_{j}\right]}, j < n \right\}, \]
and we consider the set $G_n(y,z) = \calA_n(y) \cap \calB_n(z)$. We compute the first two moments of
\[ Y_n(y,z) = \sum_{|u|=n} \ind{u \in G_n(y,z)}, \]
to bound from below $\P(Y_n(y,z) \geq 1)$, using the Cauchy-Schwarz inequality. We begin with an upper bound of the second moment of $Y_n$.
\begin{lemma}
\label{lem:secondorder}
Under assumptions \eqref{eqn:variance} and \eqref{eqn:monooptimalpath}, there exists $C>0$ such that for all $y \geq 0$, $z>0$ and $n \in \N$, we have
\[ \E(Y_n(y,z)^2) \leq C z (1 +y\indset{B}(1)) e^{-\theta y}. \]
\end{lemma}

\begin{proof}
Applying Lemma \ref{prop:spinaldecomposition}, we have
\begin{align*}
  \E(Y_n(y,z)^2)
  &= \bar{\E}\left[ \frac{1}{W_n} Y_n(y,z)^2 \right]
  = \hat{\E}\left[ \frac{1}{W_n} \sum_{|u|=n} \ind{u \in G_n(y,z)} Y_n(y,z) \right]\\
  &= \hat{\E}\left[ e^{-\theta V(w_n) + K^{(n)}_n} \ind{w_n \in G_n(y,z)} Y_n(y,z) \right].
\end{align*}
Using the fact that $w_n \in \calA_n(y) \subset G_n(y,z)$, we have
\[
  \E(Y_n(y,z)^2) \leq Cn^{\theta \lambda} e^{-\theta y} \hat{\E}\left[ Y_n(y,z) \ind{w_n \in G_n(y,z)} \right].
\]

We decompose $Y_n(y,z)$ along the spine, to obtain
\[
  Y_n(y,z) \leq \ind{w_n \in G_n(y,z)} + \sum_{k=0}^{n-1} \sum_{u \in \Omega(w_k)} Y_n(u,y),
\]
where, for $u \in \T$ and $y \geq 0$, we write $Y_n(u,y) = \sum_{|u'|=n, u'>u} \ind{u' \in \calA_n(y)}$. Let $k < n$. We recall that conditionally on $\calG_n$, the branching random walks of the descendants of distinct children $u,v \in \Omega(w_k)$ are independent. Moreover, the branching random walk starting from an individual $u \in \Omega(w_k)$ has law $\P_{V(u),k+1}$. As a consequence, for $y \geq 0$, $k < n$ and $u \in \Omega(w_k)$,
\[
  \hat{\E}\left[ Y_n(u,y)| \calG_n \right]
  = \E_{V(u),k+1}\left[ \sum_{|u'|=n-k-1} \ind{V(u') \geq m_n + y} \ind{V(u'_j) \leq g^{(n)}_{k+j+1}+y, j \leq n-k}\right].
\]
We use \eqref{eqn:manytoone} and \eqref{eqn:monoestimate} to obtain
\[
  \hat{\E}\left[ Y_n(u,y)| \calG_n \right]
  \leq C n^{\lambda \theta} e^{-\theta y} e^{\theta V(u)-K^{(n)}_{k+1}}\P_{V(u),k+1}\left(\begin{array}{l}  S_j \leq g^{(n)}_{j+k+1}+y, j\leq n-k-1 \\ S_{n-k-1} \geq m_n + y \end{array}\right).
\]
We now apply Lemma \ref{lem:ballot_general}. For all $p \leq P$ and $k \in [\alpha^{(n)}_{p-1}, \alpha^{(n)}_p)$, we have
\begin{multline}
  \P_{V(u),k+1}(S_{n-k-1} \geq m_n + y, S_j \leq g^{(n)}_{j+k+1}+y, j \leq n-k-1)\\
   \label{eqn:ending}
   \leq
  \begin{cases}
    C\frac{1+(g^{(n)}_{k+1}+y-V(u))_+\indset{B}(p)}{(\alpha^{(n)}_p-k+1)^{\indset{B}(p)/2} n^{(1 + \indset{B}(p))/2}} & \mathrm{if} \quad p < P-1\\
    C \frac{1 + (g^{(n)}_{k+1} + y - V(u))_+\indset{B}(P)}{(n-k+1)^{1/2 + \indset{B}(P)}} & \mathrm{if} \quad p = P.
  \end{cases}
\end{multline}

Let $p \leq P$ and $k \in [\alpha^{(n)}_{p-1}, \alpha^{(n)}_p)$, we compute the quantity
\[
  h_k := \hat{\E} \left[ \ind{w_n \in G_n(y,z)}\sum_{u \in \Omega(w_k)} (1 + (g^{(n)}_{k+1}+y-V(u))_+\indset{B}(p)) e^{\theta V(u) - g^{(n)}_{k+1}} \right].
\]
Using \eqref{eqn:monolegendre}, the definition of $\xi(w_{k})$ and the fact $x \mapsto x_+$ is Lipschitz, we have
\begin{align*}
  h_k
  &\leq C \hat{\E}\left[ e^{\theta (V(w_{k}) - g^{(n)}_{k})} (1 + (g^{(n)}_{k}+y-V(w_{k})_+) \xi(w_{k}) \ind{w_n \in G_n(y,z)} \right]\\
 &\leq C z \hat{\E}\left[ e^{\frac{\theta}{2}(V(w_k)-g^{(n)}_k)}(1 + (g^{(n)}_{k}+y-V(w_{k})_+) \xi(w_{k}) \ind{w_n \in \calA_n(y,z)} \right]
\end{align*}
as $w_n \in \calB_n(z)$.
Decomposing this expectation with respect tot the value taken by $V(w_k)$, we obtain
\[  h_k \leq C z e^{\theta y}
   \sum_{i=0}^{+\infty} (1+i)e^{-\theta i/2} \P\left[ \begin{array}{l}  S_j \leq g^{(n)}_j + y, j \in B^{(n)}\\ S_n \geq m_n+y, S_{k} - g^{(n)}_{k} - y \in [-i-1,-i]\end{array}\right].
\]
We apply the Markov property at time $k$ and Lemma \ref{lem:ballot_general} to obtain, if $p \in B$
\begin{equation}
  \label{eqn:beginning}
  h_k \leq
  \begin{cases}
    C z \frac{(1+y)e^{\theta y}}{k^{3/2}(\alpha^{(n)}_1-k+1)^{1/2}} \frac{n^{1/2}}{n^{\theta \lambda}} & \mathrm{if} \quad p=1\\
    C z \frac{(1+y\indset{B}(1))e^{\theta y}}{(k-\alpha^{(n)}_{p-1})^{1/2}(\alpha^{(n)}_p-k+1)^{1/2} n^{1/2}} \frac{1}{n^{\theta \lambda}} & \mathrm{if} \quad 1 < p < P\\
    C z \frac{(1 + y\indset{B}(1))e^{\theta y}}{(k-\alpha^{(n)}_{P-1}+1)^{1/2}(n-k+1)^{3/2}} \frac{1}{n^{\theta \lambda}} & \mathrm{if} \quad p = P.
  \end{cases}
\end{equation}
In the same way, if $p \not \in B$, we have
\begin{equation}
  h_k \leq 
  \begin{cases}
    C z e^{\theta y} \frac{1}{k^{1/2}} \frac{1}{n^{\theta \lambda}} & \mathrm{if} \quad k < \alpha^{(n)}_1\\
    C z e^{\theta y} \frac{1+y\indset{B}(1)}{n^{1/2}} \frac{1}{n^{\theta \lambda}} & \mathrm{if} \quad \alpha^{(n)}_1 \leq k < \alpha^{(n)}_{P-1}\\
    C z e^{\theta y} \frac{1+y\indset{B}(1)}{(n-k+1)^{1/2}} \frac{1}{n^{\theta \lambda}} & \mathrm{otherwise,}
  \end{cases}
\label{eqn:beginningprime}
\end{equation}
applying again Lemma \ref{lem:ballot_general}.

For $p \leq P$ we denote by
\[  H_p := \sum_{k = \alpha^{(n)}_{p-1}}^{\alpha^{(n)}_p-1} \hat{\E}\left[  \ind{w_n \in G_n(y,z)} \sum_{u\in \Omega(w_k)} Y_n(u,y) \right] \leq C \sum_{k = \alpha^{(n)}_{p-1}}^{\alpha^{(n)}_p-1} h_k e^{\theta g^{(n)}_{k+1} - K^{(n)}_{k+1}}.
\]
Using \eqref{eqn:monolegendre}, and summing the estimates \eqref{eqn:ending}, \eqref{eqn:beginning} and \eqref{eqn:beginningprime}, we have $H_p \leq C z (1 + y\indset{B}(1))e^{-\theta y}$ for any $p \leq P$. To conclude this proof, we observe that
\[
  \E(Y_n(y,z)^2)  \leq \sum_{p=1}^P H_p + Cn^{\theta \lambda} e^{-\theta y}\P(w_n \in G_n(y,z))
  \leq C z (1 + y\indset{B}(1))e^{-\theta y},
\]
as Lemma \ref{lem:ballot_general} implies $\P(w_n \in G_n(y,z)) \leq C (1 + y \indset{B}(1))n^{-\theta \lambda}$.
\end{proof}

We now prove the following result, a lower bound on the first moment of $Y_n(y,B)$.
\begin{lemma}
\label{lem:firstorder}
Under assumptions \eqref{eqn:variance}, \eqref{eqn:integrability} and \eqref{eqn:monooptimalpath}, there exists $c>0$ and $z>0$ such that for any $n \geq 0$ and $y \in [0, \sqrt{n}]$, we have $\E(Y_n(y,z)) \geq c (1 + \indset{B}(1)y)e^{-\theta y}$.
\end{lemma}

\begin{proof}
Using Lemma \ref{prop:spinaldecomposition}, we have
\[  \E(Y_n(y,z)) = \hat{\E}\left[ e^{-\theta V(w_n) + K^{(n)}_n} \ind{w_n \in G_n(y,z)} \right] \geq c n^{\theta \lambda}e^{-\theta y} \hat{\P}(w_n \in G_n(y,z)).\]
We observe that $\hat{\P} (w_n \in G_n(y,z)) = \hat{\P}(w_n \in \calA_n(y)) - \hat{\P} (w_n \in \calA_n(y) \cap \calB_n(z)^c)$. By Lemma \ref{lem:ballotlowerbound}, for any $n \geq 1$ and $y \in [0,\sqrt{n}]$
\[
 \hat{\P}(w_n \in \calA_n(y)) = \P(S_n \geq m_n + y, S_j \leq g^{(n)}_j+y, j \leq n)\geq c (1 + y\indset{B}(1))n^{-\theta \lambda}.
\]
Therefore, we only need to bound from above $\hat{\P} (w_n \in \calA_n(y) \cap \calB_n(z)^c)$ for $z>0$ large enough.

We denote by $\tau^{(n)}(z) = \inf\left\{k \leq n : \xi(w_k) \geq z \exp\left(-\frac{\theta}{2}\left[V(w_k)-g^{(n)}_{k}\right]\right) \right\}$, and for any $p \leq P$, we set $\pi_p =\hat{\P}\left( w_n \in \calA_n(y),\tau^{(n)}(z) \in(\alpha^{(n)}_{p-1}, \alpha^{(n)}_{p}] \right)$. We introduce the random variables
\[
  (\xi_p, \Delta_p) \egaldistr \left( V(w_{k+1}-V(w_k), \xi(w_k) \right) \quad \text{for } 
k \in [\alpha^{(n)}_{p-1}, \alpha^{(n)}_p).\]
Let $(\xi^p_n,\Delta^p_n)$ be i.i.d random variables with the same law as $(\xi_p,\Delta_p)$, we write $T^p_n = \Delta^p_1 + \cdots + \Delta^p_n$. For $p \in B$, we introduce
\[ \chi_p : z \longmapsto \hat{\E}\left[ \left( 1 + \left( \log_+(\xi_p)-\Delta_p - \log z\right)_+ \right)^2 \ind{\xi_p \geq z} \right],\]
and for $p \in B^c$,
\begin{multline*}
  \qquad \tilde{\chi}_p : z \longmapsto \hat{\E}\left[ \frac{\left(\log_+ (\xi_p)-\Delta_p - \log z/2 \right)_+}{\delta}\right]
  \\
  + \begin{cases}
    \E\left[ \sum_{k=0}^{+\infty} \ind{T^p_k \geq (\delta k + \log z) / 2} \right] & \mathrm{if} \quad E_p(\theta) < 0\\
    \E\left[ \sum_{k=0}^{+\infty} \ind{T^p_k \leq -(\delta k + \log z)/2} \right] & \mathrm{if} \quad E_p(\theta) = 0, E_{p-1}(\theta) < 0.
  \end{cases}\qquad
\end{multline*}

First, if $p=1$, we apply the Markov property at time $\alpha^{(n)}_1$ and Lemma \ref{lem:ballot_general} to obtain
\[ \pi_1 \leq C \frac{1}{n^{(1+\indset{B}(P))/2}} \P\left(T^1_j \leq g^{(n)}_j +y, j \leq \alpha^{(n)}_1, \exists k \leq \alpha^{(n)}_1 : \xi^1_k \geq ze^{-\theta/2 (T^1_k-g^{(n)}_k)}\right). \]
As a consequence, if $1 \in B$, we apply Lemma \ref{lem:ballot_spine} to obtain $\pi_1 \leq C \frac{1+y}{n^{\theta \lambda}} \chi_1(z)$ ; and if $1 \not \in B$, then $E_1<0$ so, applying Lemma \ref{lem:hsurobbins_spine} we have $\pi_1 \leq C \frac{1}{n^{\theta \lambda}} \tilde{\chi}_1(z)$.

We now suppose that $1< p < P$. Applying the Markov property at times $\alpha^{(n)}_p$ and $\alpha^{(n)}_{p-1}$, we have
\begin{equation}
 \label{eqn:twomarkov}
 \pi_p \leq  \frac{C}{n^{(1+\indset{B}(P))/2}} \hat{\E}\left[ \ind{V(w_j) \leq g^{(n)}_j + y, j \leq \alpha^{(n)}_{p-1}} \phi_p\left(V(w_{\alpha^{(n)}_{p-1}})\right) \right],
\end{equation}
where we write, for $s \in \R$
\[
  \phi_p(s) = \P_s \left[ T^p_j \leq g^{(n)}_{\alpha^{(n)}_{p-1}+j}+y, j \leq \alpha^{(n)}_p-\alpha^{(n)}_{p-1}, \tau^{(n)}(z) \in (\alpha^{(n)}_{p-1}, \alpha^{(n)}_p] \right].
\]
If $p \in B$, applying Lemma \ref{lem:ballot_spine}, we have $\phi_p(s) \leq \frac{1 + y + s}{n^{1/2}} \chi_p(z)$,
and, by Theorem \ref{thm:locallimit_excursion},
\[ \sup_{n \in \N} \frac{1}{n^{1/2}}\E\left[ \left.\left|S_{\alpha^{(n)}_{p-1}}-\bar{a}^{(n)}_{\alpha^{(n)}_{p-1}}\right| \right| S_j \leq g^{(n)}_j + y,j \leq \alpha^{(n)}_{p-1}\right] < +\infty. \]
By Lemma \ref{lem:ballot_general}, as $y \leq \sqrt{n}$, we have $\pi_p \leq \frac{C(1 + y\indset{B}(1))}{n^{\theta \lambda}} \chi_p(z)$.

In the same way, if $p \not \in B$, we use Lemma \ref{lem:hsurobbins_spine} --as well as time-reversal when $E_p(\theta)=0$ and $E_{p-1}(\theta)<0$-- to have $\phi_p(s) \leq \tilde{\chi}_p(z)$, which, thanks to \eqref{eqn:twomarkov} leads to $\pi_p \leq \frac{1+y\indset{B}(1)}{n^{\theta\lambda}} \tilde{\chi}_p(z)$.

If $P \in B$, we apply the Markov property and Lemma \ref{lem:excursion_spine} to obtain
\[
  \pi_P \leq C \E\left[\frac{1 + (S_{\alpha^{(n)}_{P-1}}-a^{(n)}_{\alpha^{(n)}_{P-1}} + y)_+}{n^{3/2}} \ind{S_j \leq g^{(n)}_j+y , j \leq \alpha^{(n)}_{P-1}} \right] \leq C \frac{1 + y\indset{B}(1)}{n^{\theta \lambda}} \chi_P(z).
\]
If $P \not \in B$, we use the time-reversal, then Lemma \ref{lem:hsurobbins_spine} to obtain
\[
  \pi_P \leq C \tilde{\chi}_P(z) \sup_{h \in \R} \P\left[ S_{\alpha^{(n)}_{p-1}} \in [h,h+1], S_j \leq g^{(n)}_j +y , j \leq \alpha^{(n)}_{P-1} \right] \leq C \frac{1 + y \indset{B}(1)}{n^{\theta \lambda}} \tilde{\chi}_P(z).
\]

We conclude there exists $C>0$ such that
\[ \hat{\P}(w_n \in \calA^{(n)}(y) \cap \calB^{(n)}(z)^c) \leq C \frac{1 + y \indset{B}(1)}{n^{\theta \lambda}} \left[\sum_{p \in B} \chi_p(z) + \sum_{p \in B^c} \tilde{\chi}_p(z)\right]. \]
If $p \in B$, by \eqref{eqn:integrability} and \eqref{eqn:variance}, $\E((\log \xi_p -\Delta_p)^2)<+\infty$. In the same way, if $p \not \in B$, using \eqref{eqn:integrability} and \eqref{eqn:variance} again, we have  $\E((\log \xi_p -\Delta_p)_+)<+\infty$. Applying the dominated convergence theorem, we have  $\lim_{z \to + \infty} \sum_{p \in B} \chi_p(z) + \sum_{p \in B^c} \chi_p(z) = 0$.
Consequently, there exists $z \geq 0$ large enough such that $\hat{\P}(w_n \in \calA^{(n)}(y) \cap \calB^{(n)}(z)^c) \leq c/2(1 + y \indset{B}(1))n^{-\theta \lambda}$. Therefore
\begin{align*}
  \hat{\P}(w_n \in \calA^{(n)}(y) \cap \calB^{(n)}(z))
  &\geq \hat{\P}(w_n \in \calA^{(n)}(y)) - \hat{\P}(w_n \in \calA^{(n)}(y) \cap \calB^{(n)}(z)^c) \\
  &\geq c (1 + y \indset{B}(1))n^{-\theta \lambda}/2,
\end{align*}
which ends the proof.
\end{proof}

Using these two lemmas, we obtain a lower bound on $M_n$.
\begin{proof}[Lower bound in Theorem \ref{thm:light}]
By Lemma \ref{lem:firstorder}, there exist $c>0$ and $z>0$ such that for any $n \geq 1$ and $y \in [0,\sqrt{n}]$, we have $\E(Y_n(y,z)) \geq c (1 + y\indset{B}(1))e^{-\theta y}$. Thus, using Lemma \ref{lem:secondorder} and the Cauchy-Schwarz inequality, we have
\begin{align*}
  \P(Y_n(y,z) \geq 1) &\geq \frac{\E(Y_n(y,z))^2}{\E(Y_n(y,z)^2)} \geq \frac{\left(c(1 + y\indset{B}(1))e^{-\theta y}\right)^2}{C z (1 + y\indset{B}(1)) e^{-\theta y}}\\
  &\geq c (1 +y \indset{B}(1))e^{-\theta y}.
\end{align*}
\end{proof}

\subsection{Extension to the multi-parameter branching random walk}
\label{subsec:extension}

In this section, we extend Theorem \ref{thm:light} to BRWis such that $\bbtheta$ is non-constant, reasoning by induction on the number $T$ of different values taken by the sequence.

\begin{proof}[Proof of Theorem \ref{thm:tailestimate}]
We observe first that if $T=1$, then the branching random walk satisfies all the hypotheses of Theorem \ref{thm:light}, with optimal path $\bba$, and parameter $\theta=\phi_1$, by Proposition \ref{prop:optimization_problem}. The initiation of the recurrence is then given by Theorem \ref{thm:light}. Therefore, we only need to prove the induction hypothesis.

Let $T \in \N$, we assume that for all BRWis such that $\#\{\theta_p, p \leq P\}<T$, Theorem \ref{thm:tailestimate} holds. For $n \in \N$, we now consider a BRWis $(\T^{(n)}, V^{(n)})$ of length $n$. We write $\bba$ the optimal solution of Proposition \ref{prop:optimization_problem}, and $\theta_p = \kappa'_p(a_p)$. We assume that $T=\#\{\theta_p, p \leq P\}$, and write $\phi_1 < \phi_2 < \cdots < \phi_T$ these values, listed in the increasing order. For any $t \leq T$, let $f_t = \min\{p \leq P : \theta_p = \phi_t\}$ and $l_t = \max\{p \leq P : \theta_p = \phi_t\}$. Finally, we write $v_\mathrm{is}$ and $\lambda$ the speed and correction as defined in \eqref{eqn:speedisdef} and \eqref{eqn:logisdef}, and $m_n = nv_\mathrm{is} - \lambda \log n$ the expected position of the maximal displacement $M_n$. We now divide this BRWis into two parts, before and after the first time $\alpha_{l_1}$ such that $\theta_{{l_1}+1} > \theta_{l_1}$.

We write $l=l_1$, $v_1 = \sum_{p=1}^l (\alpha_p - \alpha_{p-1}) a_p$ and $\lambda_1 = \frac{1}{2\phi_1}\left( 1 + \indset{B}(1) + \indset{B}(l) \right)$. We denote by $\T^{(n)}_1 = \{ u \in \T : |u| \leq \alpha_l n \}$ the tree cut at generation $n_1 = \floor{\alpha_l n}$. By Proposition \ref{prop:optimization_problem}, we observe that $(\T^{(n)}_1, V^{(n)}_{|\T^{(n)}_1})$ is a BRWis which satisfies the hypotheses of Theorem \ref{thm:light}, with parameter $\theta := \phi_1$. Therefore, if we write $m^1_n = v_1 n - \lambda_1 \log n$ and $M^1_n = \max_{|u|=n_1} V(u)$, there exist $c,C>0$ such that for all $n \in \N$ large enough and $y \in [0,n^{1/2}]$, we have
\[
  c(1 + y\ind{\kappa^*_1(a_1)=0}) e^{-\phi_1 y} \leq \P(M^1_n \geq m^1_n + y) \leq C (1 + y\ind{\kappa^*_1(a_1)=0})e^{-\phi_1 y}.
\]

We now consider a branching random walk $(T^{(n)}_\text{tail}, V^{(n)}_\text{tail})$ of law $\P_{n_1,0}$, which has the law of the branching random walk of the descendants of any individual alive at time $n_1$. We write $n_\text{tail}=n - n_1$ the length of this BRWis, $v_\text{tail} = v - v_1$, $\lambda_\text{tail} = \lambda - \lambda_1$ and $m_n^\text{tail} = v_\text{tail} n - \lambda_\text{tail} \log n$.
We observe easily, by Proposition \ref{prop:optimization_problem} again, that this marked tree is a BRWis, and its optimal path is the path driven by $(a_{l+1}, \ldots, a_P)$. Moreover, $\#\{\theta_p, l<p\leq P\}=T-1<T$. Therefore, by the induction hypothesis, writing $M^\text{tail}_n = \max_{|u|=n_\text{tail}} V^\text{tail}(u)$, there exist $c,C>0$ such that for all $n \in \N$ large enough and $y \in [0,n^{1/2}]$, we have
\[
  c e^{-\phi_2 y} \leq \P(M^\text{tail}_n \geq m^\text{tail}_n + y) \leq C (1 + y)e^{-\phi_2 y}.
\]

To obtain the lower bound of Theorem \ref{thm:tailestimate}, we observe that if $M^1_n \geq m^1_n + y$, and if one of the descendants of the rightmost individual at time $n_1$ makes a displacement greater than $m^\text{tail}_n$, then $M_n \geq m_n + y$. Therefore
\[
  \P(M_n \geq nv_\mathrm{is} - \lambda \log n + y) \geq c (1 + y\indset{B}(1))e^{-\phi_1 l}
\]
for $n \in \N$ large enough and $y \in [0,n^{1/2}]$. To obtain an upper bound for $\P(M_n \geq m_n + y)$, we decompose the $n^\mathrm{th}$ generation of the branching random walk with respect to the position of their ancestors alive at time $n_1$. We write
\[ X^{(n)}(y,h) = \sum_{|u| = n_1} \ind{V(u_j) \leq f^{(n_1)}_j + y, j \in C^{(n_1)}_{1}} \ind{V(u) -m^1_n - y \in [-h-1,-h]}, \]
and, by union bound and the Markov property, we have
\begin{align*}
  &\P(M_n \geq m_n + y)\\
  \leq &\P(\exists |u| \in C^{(n_1)}_1 : V(x) \geq f^{(n_1)}_k+r^{(n)}_k + y)
   + \sum_{h =0}^{+\infty} \E(X^{(n)}(y,h)) \P(M^\text{tail}_{n} \geq m^\text{tail}_n + h).
\end{align*}
As a consequence, applying Lemma \ref{lem:upperestimate} and the upper bound of Theorem \ref{thm:light},
\begin{align*}
 \P(M_n \geq m_n + y)
  &\leq C (1 + y \ind{\kappa^*_1(a_1)=0}) e^{-\phi_1 y}\left[ 1 + \sum_{h=0}^{+\infty} (1+h)e^{(\phi_1-\phi_2) h}\right]\\
  &\leq C (1  + y \ind{\kappa^*_1(a_1)=0}) e^{-\phi_1 y},
\end{align*}
which gives the correct upper bound.
\end{proof}

\subsection{Proof of Theorem \ref{thm:main}}
\label{subsec:proof}

Using Theorem \ref{thm:tailestimate}, we are able to obtain Theorem \ref{thm:main}. To do so, we need to strengthen the estimate $\P(M_n \geq nv_\mathrm{is} - \lambda \log n)>c>0$ in something like
\[
  \lim_{y \to -\infty} \liminf_{n \to +\infty} \P(M_n \geq nv_\mathrm{is} - \lambda \log n) = 1.
\]
To do so, we will use a standard cutting argument. We use the fact that with high probability, there will be a large number of individuals alive at a fixed generation $k$, each of which having positive probability to make a descendant at generation $n$ with a displacement greater than $m_n$. Using the law of large numbers, this will be enough to conclude.

\begin{proof}[Proof of Theorem \ref{thm:main}]
Let $(\T,V)$ be a BRWis of length $n$, satisfying all hypotheses of Theorem \ref{thm:main}. To prove that the sequence $(M_n - m_n)$ is tight, we need to prove that
\[
  \lim_{K \to +\infty} \sup_{n \in \N} \P( |M_n - m_n| \geq K) = 0.
\]
By Theorem \ref{thm:tailestimate}, there exists $C>0$ such that
\[
  \sup_{n \in \N} \P(M_n \geq m_n + K) \leq C (1+K)e^{-\phi_1 K},
\]
therefore the upper bound is easy to obtain.

We now turn to the lower bound. Applying Theorem \ref{thm:tailestimate}, there exists $c_1>0$ such that
\[
  \inf_{n \in \N} \P(M_n \geq m_n) \geq c_1.
\]

Let $L_1$ be a point process of law $\calL_1$. By \eqref{eqn:breeding}, there exists $h>0$ and $N \in \N$ such that
\[
  m = \E\left[ \max\left( N, \sum_{\ell \in L_1} \ind{\ell \geq -h} \right) \right]>1.
\]
We write $\mu$ the law of $\max\left( N, \sum_{\ell \in L_1} \ind{\ell \geq -h} \right)$, and $(Z_n,n \geq 0)$ a Galton-Watson process with reproduction law $\mu$. We can easily couple $(Z_n)$ and a branching random walk $(\T_1,V_1)$ with reproduction law $\calL_1$ in such a way that for all $n \in \N$, $\sum_{|u|=n} \ind{V_1(u) \geq - n h} \geq Z_n$. By standard Galton-Watson processes theory, there exists $c_2>0$ and $\delta > 0$ such that $\inf_{k \in \N} \P(Z_k \geq \delta m^k) > c_2$.

Let $\epsilon>0$ and $R>0$ be such that $(1-c_1)^R\leq \epsilon$. We now choose $k \in \N$ such that $\delta m^k \geq R$. For any $n \in \N$, we write $u_n = (1,\ldots 1) \in \calU$. By \eqref{eqn:breeding}, for all $n \in \N$ we have $u_n \in \T$. We write $\tau$ the first time $n$ such that $u_n$ has a sibling at distance smaller than $h$, and this child has at least $R$ descendants alive at time $n+k$ whose relative position is less that $-kh$. According to the previous computations, $\tau$ is stochastically dominated by a Geometric random variable. Therefore, it exists $\tau_0 \in \N$ such that $\P(\tau > \tau_0)< \epsilon$.

Therefore, with probability at least $1-2\epsilon$, there are at least $R$ individuals alive at some time before $\tau_0 + k$, all of which are above $\inf_{j \leq \tau_0} V(u_j) - k h$. Each of these individuals $u$ starts an independent BRWis with law $\P_{k,V(u)}$, thus, using Theorem \ref{thm:tailestimate}, there exists $y>0$ such that, for all $n \geq 1$ large enough
\[
  \P(M_{n + \frac{k+\tau_0}{\alpha_1}} \geq m_n - y) \geq 1 - 4 \epsilon
\]
which ends the proof of the lower bound.
\end{proof}

\section{Phase-transition in the branching random walk with one interface}
\label{sec:twoenvi}
We consider in this section a BRWis with a single interface $(\T,V)$, or in other words, such that $P=2$. The process studied by Fang and Zeitouni in \cite{FaZ12a} can be described this way.Let $\calL_1$ and $\calL_2$ be two point processes, verifying \eqref{eqn:breeding}, and $\alpha_1 \in (0,1)$. We assume \eqref{eqn:regularity}, i.e. there exist $\bar{\theta}_1, \bar{\theta}_2$ such that for all $i \in \{1,2\}$,
\[
  \bar{\theta}_i \kappa'_i(\bar{\theta }_i) - \kappa_i(\bar{\theta}_i)=0.
\]
We also suppose there exists $\theta>0$ such that $\kappa_1$ and $\kappa_2$ are differentiable at point $\theta$ and
\begin{equation}
  \label{eqn:twoenvitheta}
  \theta(\alpha_1 \kappa'_1(\theta) + (1-\alpha_1) \kappa'_2(\theta)) - (\alpha_1 \kappa_1(\theta) + (1-\alpha_1)\kappa_2(\theta)) = 0.
\end{equation}

For all $i \in \{1,2\}$, $\kappa_i$ is a convex function on $\{\theta > 0 : \kappa_i(\theta)<+\infty\}$, which is twice differentiable on the interior of this set. As a consequence, $\theta \kappa'_i(\theta) - \kappa_i(\theta)$ is a decreasing function. Thus, $\theta$ is always between $\bar{\theta}_1$ and $\bar{\theta}_2$. We write 
\begin{equation}
  \label{eqn:twospeeds}
  v_\text{fast} = \alpha_1 \kappa'_1(\theta) + (1-\alpha_1) \kappa'_2(\theta) \text{  and  }v_\text{slow} = \alpha_1 \kappa'_1(\bar{\theta}_1) + (1-\alpha_1) \kappa'_2(\bar{\theta}_2).
\end{equation}
Note that $v_\text{slow}$ is the sum of the speeds of a branching random walk with reproduction $\calL_1$ of length $n\alpha_1$, with one with reproduction $\calL_2$ of length $n\alpha_2$.

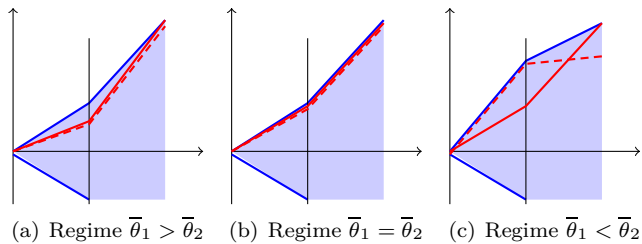
\begin{figure}[!ht]
\centering
\subfigure[Regime $\bar{\theta}_1 > \bar{\theta}_2$]
{
\label{subfig:fast}
\begin{tikzpicture}
  \fill [color=blue!20] (0,0) -- (1,0.64)--(2,1.74)--(2,-0.64) -- (1,-0.64) -- cycle ;
  \draw [->] (0,-0.7) -- (0,1.9) ;
  \draw [->] (0,0) -- (2.5,0) ;
  \draw [thick, color=blue] (0,-0.04) -- (1,-0.64);
  \draw [thick, color=blue] (0,0) -- (1,0.64)--(2,1.74);
  \draw (1,-0.7)--(1,1.5);
  \draw [thick, color=red] (0,0) -- (1,0.4)--(2,1.74);
  \draw [thick, densely dashed,color=red] (0,0) -- (1,0.36)--(2,1.66);
\end{tikzpicture}
        }
\subfigure[Regime $\bar{\theta}_1 = \bar{\theta}_2$]
{
\label{subfig:mean}
\begin{tikzpicture}
  \fill [color=blue!20] (0,0) -- (1,0.64)--(2,1.74)--(2,-0.64) -- (1,-0.64) -- cycle ;
  \draw [->] (0,-0.7) -- (0,1.9) ;
  \draw [->] (0,0) -- (2.5,0) ;
  \draw [thick, color=blue] (0,-0.04) -- (1,-0.64);
  \draw [thick, color=blue] (0,0) -- (1,0.64)--(2,1.74);
  \draw (1,-0.7)--(1,1.5);
  \draw [thick, color=red] (0,0) -- (1,0.6)--(2,1.7);
  \draw [thick, densely dashed,color=red] (0,0) -- (1,0.56)--(2,1.66);
\end{tikzpicture}
        }
\subfigure[Regime $\bar{\theta}_1 < \bar{\theta}_2$]
{
\label{subfig:slow}
\begin{tikzpicture}
  \fill [color=blue!20] (0,0) -- (1,1.2)--(2,1.7)--(2,-0.64) -- (1,-0.64) -- cycle ;
  \draw [->] (0,-0.7) -- (0,1.9) ;
  \draw [->] (0,0) -- (2.5,0) ;
  \draw [thick, color=blue] (0,-0.04) -- (1,-0.64);
  \draw [thick, color=blue] (0,0) -- (1,1.2)--(2,1.7);
  \draw (1,-0.7)--(1,1.5);
  \draw [thick, color=red] (0,0) -- (1,0.6)--(2,1.7);
  \draw [thick, densely dashed,color=red] (0,-0.04) -- (1,1.16)--(2,1.26);
\end{tikzpicture}
}
\caption{Regimes in the branching random walk with interface.}
\end{figure}

Applying Theorem \ref{thm:main} and using Proposition \ref{prop:optimization_problem}, we observe that, under \eqref{eqn:variance} and \eqref{eqn:integrability}, one of the following alternative is true.
\begin{itemize}
  \item If $\bar{\theta}_1 > \bar{\theta}_2$, then $\theta \in (\bar{\theta}_2,\bar{\theta}_1)$, $v_\text{slow} < v_\text{fast}$ and
  \[
    M_n = n v_\text{fast} - \frac{1}{2\theta} \log n + O_\P(1),
  \]
  in which case the optimal path is at time $\alpha_1 n$ at distance $O(n)$ from the frontier of the branching random walk. The rightmost individual at time $n$ is at distance $O(n)$ from the rightmost child of the rightmost individual alive at time $\alpha_1 n$ (case \ref{subfig:fast}).
  \item If $\bar{\theta_1} = \bar{\theta}_2$, then $\theta = \bar{\theta}_1 = \bar{\theta}_2$, $v_\text{slow} = v_\text{fast}$ and
  \[
    M_n = n v_\text{fast} - \frac{3}{2\theta} \log n + O_\P(1),
  \]
  and the process behaves similarly to time-homogeneous branching random walk, the path leading to the rightmost individual at time $n$ stays at any time within distance $O(\sqrt{n})$ from the frontier of the branching random walk (case \ref{subfig:mean}).
  \item If $\bar{\theta}_1 < \bar{\theta}_2$, then $v_\text{slow} < v_\text{fast}$ and
  \[
    M_n = n v_\text{slow} - \left[ \frac{3}{2\bar{\theta}_1} + \frac{3}{2\bar{\theta}_2} \right] \log n + O_\P(1),
  \]
  in other words, the logarithmic corrections add up, and the rightmost individual at time $n$ descend from one of the rightmost individuals alive at time $\alpha_1 n$ (case \ref{subfig:slow}).
\end{itemize}

\begin{figure}[!ht]
\centering
\begin{tikzpicture}
  \fill [fill=gray!20] (0,0.5)--(2.5,0.5) -- plot [smooth,domain=2.5:0.316,samples=200] (\x-0.25,1.25/\x) -- (0,3.95) -- cycle;
  \draw [thick] (0,0.5)--(2.25,0.5) -- plot [smooth,domain=2.5:0.316,samples=200] (\x-0.25,1.25/\x);
  \draw [->] (-0.5,0) -- (8,0) node [above] {$\bar{\theta}_1$};
  \draw [->] (0,-0.5) -- (0,4) node[left] {$\lambda$};
  \draw (2.25,-0.1) -- (2.25,0.1);
  \draw (2.25,-0.1)  node [below] {$\bar{\theta}_2$};
  \draw [dashed] (2.25,0) -- (2.25,1.5) node {$\bullet$};
  \draw [dashed] (0,1.5)--(8,1.5);
  \draw [thick] plot [smooth, domain=2.25:8, samples=200] (\x, 1.2+ 3.75/\x);
  \draw [decorate,decoration={brace}]
(2.25,-0.1) -- (0,-0.1) node[below=0.2cm,pos=0.5] {$v_\text{fast}$};  
  \draw [decorate,decoration={brace}]
(8,-0.1) -- (2.25,-0.1) node[below=0.2cm,pos=0.5] {$v_\text{slow}$};  
\end{tikzpicture}
\caption{Black and grey areas are the set of possible values for the logarithmic correction $\lambda$, given $\bar{\theta}_1$ and $\bar{\theta}_2$; as $\bar{\theta}_1$ grows bigger than $\bar{\theta}_2$, logarithmic correction exhibit a sharp phase transition.}
\label{fig:logcorrection}
\end{figure}
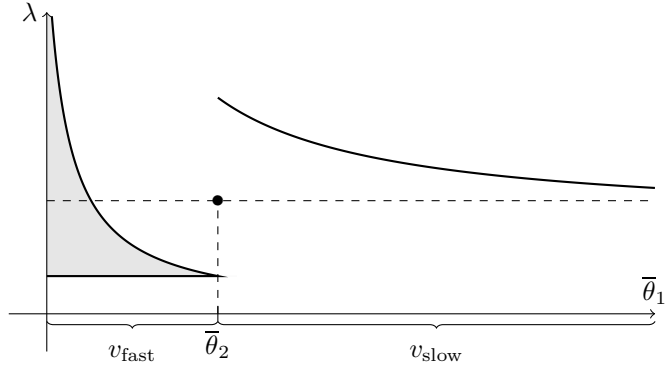

We observe, using Lagrange theorem --see Appendix \ref{app:optimization_problem}-- that
\[
v_\text{fast} = \sup\left\{ \alpha_1 a_1 + (1 - \alpha_1) a_2 : \alpha_1 \kappa^*_1(a_1) + (1-\alpha_1) \kappa^*_2(a_2) \leq 0 \right\},
\]
\[
v_\text{slow} = \sup\left\{ \alpha_1 a_1 + (1 - \alpha_1) a_2 : \alpha_1 \kappa^*_1(a_1) \leq 0, \alpha_1 \kappa^*_1(a_1) + (1-\alpha_1) \kappa^*_2(a_2) \leq 0 \right\}.
\]
Therefore, a branching random walk goes at speed $v_\text{slow}$ if the condition $\kappa^*_1(a_1) \leq 0$ matters to solve \eqref{eqn:optimization_problem}. If this is the case, the ``theoretical optimal path'' would cross the frontier of the branching random walk, thus no individual could follow it. But under these circumstances, the closer the individual is to the frontier at time $\alpha_1 n$, the better the probability that they are the ancestors of the rightmost individual at time $n$. Otherwise, at time $\alpha_1 n$, there is a large number of individuals around $\alpha_1 a_1 n$, each of which having small probability to be the rightmost individual, thus the logarithmic correction is the same as the one obtained computing the maximal displacement of a large number of independent random walks.

Although the speed $v_\mathrm{is}$ varies continuously as $\bar{\theta}_1$ grows bigger than $\bar{\theta}_2$, the logarithmic correction $\lambda$ exhibits a phase transition. We represent in Figure \ref{fig:logcorrection} the set of possible values taken by $\lambda$ for different values of $\bar{\theta}_1$ and $\bar{\theta}_2$. The frontier of this set does not depend on the value of $\alpha_1$, the position of the interface.

\appendix

\section{Time-inhomogeneous random walk estimates}
\label{app:estimates}
In this section, we prove the random walk estimates we defined in Section \ref{sec:randomwalk}.

\subsection{Proof of Lemmas \ref{lem:ballot} and \ref{lem:excursion}}
\label{app:estimates1}

We recall that $T$ is a centred random walk with finite variance. We first prove Lemma \ref{lem:ballot}: there exists $C>0$ such that $\P(T_j \geq -y -A j^\alpha, j \leq n) \leq C (1 + y)n^{-1/2}$.

\begin{proof}[Proof of Lemma \ref{lem:ballot}]
Let $\alpha \in [0,1/2)$, $A>0$ and $(f_n) \in \R^\N$ such that for all $n \in \N$, $|f_n| \leq A n^\alpha$. For all $y \geq 0$, we have
\[
  \P(T_j \geq - y - f_j, j \leq n) \leq \P( T_j \geq -y - A j^\alpha, j \leq n),
\]
thus we now bound this later probability. More precisely, for $y \geq 0$, we denote by
\[
  \tau_y = \inf\{ n \geq 0 : T_n \leq -y - A j^\alpha \},
\]
and we take interest in $\P(\tau_y \geq n)$.

For $a \in \R\backslash\{0\}$, we set
\[
  H_a =
  \begin{cases}
    \inf\{ n \geq 0 : T_n \geq a\} &\text{ if } a >0\\
    \inf\{ n \geq 0 : T_n \leq a\} &\text{ if } a < 0.
  \end{cases}
\]
By Theorem \ref{thm:ballot1}, there exists $K>0$ such that for all $a \in \R$, $\P(H_a \geq n) \leq K (1 + |a|) n^{-1/2}$. As a result, for all $y \geq n^\alpha$, we have
\begin{equation}
 \label{eqn:longshot}
  \P(\tau_y \geq n) \leq \P(H_{-(A+1)n^\alpha} \geq n) \leq K ((A+1) n^\alpha + 1) n^{-1/2} \leq K(A+1) (1 + y)n^{-1/2},
\end{equation}
hence it is enough to consider the case $y \leq n^\alpha$.

Let $\gamma \in (2\alpha,1)$, we observe that the following decomposition holds
\begin{equation}
  \label{eqn:decompositionRandomWalk}
  \P(\tau_y \geq n) \leq \P(H_{n^\alpha} \leq \min(\tau_y,n^\gamma), \tau_y \geq n) + \P(n^\gamma < \min(H_{n^\alpha},\tau_y)),
\end{equation}
and we bound these two parts separately.

We first observe that for all $y \leq n^{\alpha}$, we have
\[
  \P(n^\gamma < \min(H_{n^\alpha},\tau_y)) \leq \P(n^\gamma < \min(H_{n^\alpha},H_{-(A+1)n^\alpha})) \leq \P( \max_{j \leq n^\gamma} |T_j| \leq (A+1) n^\alpha).
\]
As a result, using Mogul'ski\u\i{} small deviations estimates \cite{Mog}, we have
\[
  \limsup_{n \to +\infty} n^{2\alpha - \gamma} \log \sup_{y \leq n^\alpha} \P(n^\gamma < \min(H_{n^\alpha},\tau_y)) \leq - \frac{\pi^2 \sigma^2}{8(A+1)^2} < 0,
\]
in particular there exists $K_1 >0$ such that for all $y \leq n^\alpha$,
\begin{equation}
  \label{eqn:decompositionSmallDeviations}
  \P(n^\gamma < \min(H_{n^\alpha},\tau_y)) \leq K_1(1+y)n^{-1/2}.
\end{equation}

We now bound the other term of Equation \ref{eqn:decompositionRandomWalk}. Applying the Markov property at time $H_{n^\alpha}$, for all $n$ large enough, we have
\begin{align*}
  \P(H_{n^\alpha} \leq \min(\tau_y,n^\gamma), \tau_y \geq n) &\leq \left(\ind{H_{n^\alpha}<\min(\tau_y,n^\gamma)} \P_{T_{H_{n^\alpha}}}(\tau_y \geq n-n^\gamma)\right)\\
  &\leq \E\left(\ind{H_{n^\alpha}<\tau_y} \P_{T_{H_{n^\alpha}}} (H_{-(A+1) n^\alpha} > n/2 ) \right)\\
  &\leq \E\left( \ind{H_{n^\alpha}<\tau_y} K \frac{(1 + (A+1)n^{\alpha} + T_{H_{n^\alpha}})}{(n/2)^{1/2}}  \right)\\
  &\leq 2 K(A+3) n^{-1/2} \E\left( \ind{H_{n^\alpha}<\tau_y} T_{H_{n^\alpha}} \right)
\end{align*}
using Theorem \ref{thm:ballot1} and the fact that $T_{H_{n^\alpha}}\geq n^\alpha \geq 1$ a.s. Moreover, as $(T_n)$ is a martingale, applying the optional stopping theorem we have
\[  0 = \E(T_{H_{n^\alpha}\wedge \tau_y}) = \E(T_{H_{n^\alpha}} \ind{H_{n^\alpha}<\tau_y}) + \E(T_{\tau_y} \ind{\tau_y <H_{n^\alpha}}),\]
therefore
\[
  \E(T_{H_{n^\alpha}} \ind{H_{n^\alpha}<\tau_y}) = \E\left( - T_{\tau_y} \ind{\tau_y < H_{n^\alpha}} \right) \leq \E\left( -T_{\tau_y}\right) := \psi(y).
\]
Using \cite[Theorem 7]{DSV}, $\psi(y)<+\infty$ for all $y \geq 0$. We conclude there exists $K_2>0$ such that for all $n \in \N$ and $y \in [0,n^\alpha]$,
\[
  \P(H_{n^\alpha} \leq \min(\tau_y,n^\gamma), \tau_y \geq n) \leq K_2 \psi(y) n^{-1/2}.
\]
To conclude the proof, it is enough to observe that $\psi$ increases at most at a linear rate.

Let $y,y' \geq 0$, as $j \mapsto - A j^\alpha$ is decreasing, we observe that
\begin{align*}
  \psi(y+y') &= - \E(T_{\tau_{y+y'}}) = -\E(T_{\tau_y} + T_{\tau_{y+y'}} - T_{\tau_y})\\
  &\leq -\E(T_{\tau_y}) - \E(T_{\tau_{y'}}) = \psi(y) + \psi(y').
\end{align*}
As $\psi$ is sub-additive, there exists $L>0$ such that $\psi(y) \leq (L+1)y$ for all $y \geq 0$. Coming back to \eqref{eqn:decompositionRandomWalk}, we conclude there exists $C>0$ such that for all $n \in \N$ and $y \leq n^\alpha$,
\[
  \P(\tau_y \geq n) \leq C (1+y)n^{-1/2}.
\]

\end{proof}

We now prove Lemma \ref{lem:excursion}: there exists $C>0$ such that for all $p,q \in \N$, $x,h \geq 0$  and $y \in \R$, we have
\begin{equation*}
  \P\left( \begin{array}{l} T_j \geq -x \ind{j \leq p} + y \ind{p<j \leq p+q}, j \leq p+q \\ T_{p+q} \in [y+h,y+h+1] \end{array} \right)  \leq C \frac{1+x}{p^{1/2}} \frac{1}{\max(p,q)^{1/2}} \frac{1+h}{q^{1/2}}.
\end{equation*}

\begin{proof}[Proof of Lemma \ref{lem:excursion}]
We denote by $p' = \floor{p/2}$, $q'=\floor{q/2}$ and by $p''=p-p'$, $q''=q-q'$. Applying the Markov property at time $p'$, we have
\begin{multline*}
  \P\left(T_{p+q}-y-h \in [0,1], T_j \geq -x \ind{j \leq p} + y \ind{p< j \leq p+q}, j \leq p+q\right)\\
  \leq \P\left(T_j \geq -x, j \leq p'\right) \sup_{z \geq -x} \P\left(T_{p''+q}+z -y-h\in [0,1], T_j+z \geq y, p'' < j \leq p''+q\right).
\end{multline*}
We set $\hat{T}_k=T_{p''+q}-T_{p''+q-k}$, which , once again, has same law as $T$. For all $z \in \R$, we have
\begin{multline*}
  \P\left(T_{p''+q} -z-h \in [0,1], T_j \geq z, p'' < j \leq p''+q\right)\\
  \leq \P\left(\hat{T}_{p''+q} \in [z+h,z+h+1], \hat{T}_j \leq h+1, j \leq q\right).
\end{multline*}
Applying again the Markov property at time $q'$, we deduce that
\begin{multline*}
  \P\left(T_{p+q} \in [y+h,y+h+1], T_j \geq -x \ind{j \leq p} + y \ind{p< j \leq p+q}, j \leq p+q\right)\\
  \leq \underbrace{\P\left(T_j \geq -x, j \leq p'\right)}_{\frac{1+x}{p^{1/2}}} \times \underbrace{\P\left(T_j \leq h+1, j \leq q'\right)}_{\frac{1+h}{q^{1/2}}} \times \sup_{z \in \R} \underbrace{\P\left(T_{p''+q''} \in [z,z+1]\right)}_{\frac{1}{\max(p,q)^{1/2}}},
\end{multline*}
using Theorems \ref{thm:locallimit} and \ref{thm:ballot1}.
\end{proof}

\subsection{Proof of Lemma \ref{lem:ballot_general}}
\label{app:estimates2}

We recall here the notations of Lemma \ref{lem:ballot_general}. Let $p,q,r \in \N$, set $n=p+q+r$. The time-inhomogeneous random walk $S$ consists of $p$ steps of independent centred random walk with finite variance, $q$ steps of independent random variables, then $r$ steps of another centred random walk with finite variance.  Let $A \in \R$, and $x,y \in \R_+$, $h \in \R$, we write
\begin{align*}
  \Gamma^{A,1}(x,y,h) &= \left\{ s \in \R^n : \forall k \leq p, s_k \geq - x \right\}\\
  \Gamma^{A,3}(x,y,h) &= \left\{ s \in \R^n : \forall k \in [n-r,n], s_k \geq y + A \log \tfrac{n}{n-k+1} \right\}.
\end{align*}

\begin{proof}[Proof of Lemma \ref{lem:ballot_general}]
Let $A>0$, $p,q,r \in \N$, $y \geq 0$ and $h \in \R$. Without loss of generality, we can assume that both $p$ and $r$ are even (by changing $q$ in $q+1$ or $q+2$). 

If $F=\emptyset$, Lemma \ref{lem:ballot_general} is an easy consequence of Theorem \ref{thm:locallimit}.

If $F = \{1\}$, applying the Markov property at time $p/2$, we obtain
\begin{multline*}
  \P\left[S_n + A \log n \in [y+h,y+h+1], (S_k, k \leq n) \in \Gamma^{A,1}(x,y,h)\right]\\
   \leq \P\left(S_j \geq -x, j \leq p/2\right) \sup_{z \in \R} \P\left(S_n-S_{p/2} \in [z,z+1]\right)\leq C \frac{1+x}{p^{1/2}} \frac{1}{\max(p,r)^{1/2}},
\end{multline*}
using Theorems \ref{thm:ballot} and \ref{thm:locallimit} respectivelly. If $F=\{3\}$, we apply the time-reversal, let $\hat{S}_j = S_n - S_{n-j}$. We have
\begin{multline*}
  \P\left[ S_n + A \log n-y-h \in [0,1], S_j \geq  y + A \log \tfrac{n}{n-j+1}, n-r \leq j \leq n \right]\\
  \leq \P\left[ \hat{S}_n + A \log n -y-h\in [0,1], \hat{S}_j \leq h + 1 - A \log (j+1), j \leq r \right]\\
  \leq C \frac{1 + h_+}{r^{1/2}}\frac{1}{\max(p,r)^{1/2}},
\end{multline*}
by the same arguments as above.

Finally, if $F=\{1,3\}$, applying Markov property at time $p/2$, and time-reversal
\begin{multline*}
 \P\left[ (S_k, k \leq n) \in \Gamma^{A,1}(x,y,h) \cap \Gamma^{A,3}(x,y,h)\right]\\
  \leq\P\left[ S_j \geq -x, j \leq p/2 \right] \sup_{z  \in \R} \P\left[ \hat{S}_{n-p/2}\in [z,z+1], \hat{S}_j \leq h +1 - A \log (j+1), j \leq r \right].
\end{multline*}
As a consequence, using once again the same arguments
\[ \P\left[ (S_k, k \leq n) \in \Gamma^{A,1}(x,y,h) \cap \Gamma^{A,3}(x,y,h) \right] \leq C\frac{1+x}{p^{1/2}} \frac{1}{\max(p,r)^{1/2}} \frac{1+h_+}{r^{1/2}}.\] 
\end{proof}

\subsection{Proof of Lemma \ref{lem:ballotlowerbound}}
\label{app:estimates3}

We consider a collection of independent random variables $(X_n^p,n \geq 0, p \leq P)$, with, for all $p \leq P$, $(X_n^p,n \geq 0)$ an i.i.d. sequence of real-valued centred random variables with finite variance. Let $n \geq 1$, we write, for $k \leq n$, $S_k = \sum_{p=1}^P \sum_{j=1}^k X_j\ind{j \in (\alpha^{(n)}_{p-1},\alpha^{(n)}_p]}$. For $F \subset \{1, 3\}$ and $x,y,\delta \in \R_+$, we write
\[
 \Upsilon^F(x,y,\delta) = \left\{ s \in \R^n :
 \begin{array}{l}
   s_k \geq - x \ind{1 \in F} - \delta k \ind{1 \not \in F},  k \leq \alpha^{(n)}_{1}, s_k \geq 0, k\in (\alpha^{(n)}_1,\alpha^{(n)}_{P-1}]\\
   s_k \geq y\ind{3 \in F}-\delta (n-k) \ind{3 \not \in F},k \in (\alpha^{(n)}_{p-1},n]
 \end{array}
 \right\}.
\]
There exists $c>0$ such that for any $F \subset \{1,3\}$, $x \in [0,n^{1/2}]$, $y \in [-n^{1/2},n^{1/2}]$ and $\delta > 0$
\[ \P\left[ S_n \leq y+1, S \in \Upsilon^F(x,y,\delta)\right] \geq c \frac{1 + x\indset{F}(1)}{n^{\indset{F}(1)/2}} \frac{1}{n^{1/2}} \frac{1}{n^{\indset{F}(3)/2}}. \]

\begin{proof}[Proof of Lemma \ref{lem:ballotlowerbound}]
Let $n \geq 1$, $x, |y| \in [0,n^{1/2}]$ and $\delta > 0$. We denote by
\[ \Omega^F(\delta, y) = \left\{ s \in \R^{n-\alpha^{(n)}_1} :  \begin{array}{l}
   \forall k\leq \alpha^{(n)}_{P-1}-\alpha^{(n)}_1,  s_k \geq 0\\
   \forall k \in (\alpha^{(n)}_{p-1},n],  s_k \geq y\ind{3 \in F}-\delta (n-k) \ind{3 \not \in F}
 \end{array}
 \right\}.
 \]
Applying the Markov property at time $\alpha^{(n)}_1$, we have
\begin{multline*}
  \P\left[ S_n \leq y+1, S \in \Upsilon^F(x,y,\delta) \right]\\
  = \E\left[ \ind{S_j \geq -x \ind{1 \in F} - \delta k \ind{1 \not \in F}} \P_{\alpha^{(n)}_1, S{\alpha^{(n)}_1}} \left( S_{n-\alpha^{(n)}_1} \leq y+1, S \in \Omega^F(\delta, y) \right)  \right].
\end{multline*}

On the one hand, if $1 \in F$, we have
\begin{multline*}
  \P\left[ S_n \leq y+1, S \in \Upsilon^F(x,y,\delta) \right] \\
  \geq \P\left(S_j \geq -x, S_{\alpha^{(n)}_1} \in [3n^{1/2}, 4n^{1/2}]\right)  \qquad \qquad \qquad \qquad \qquad \qquad \\
  \times\inf_{u \in [3n^{1/2},4n^{1/2}]}\P_{\alpha^{(n)}_1, u} \left( S_{n-\alpha^{(n)}_1} \leq y+1, S \in \Omega^F(\delta, y) \right).
\end{multline*}
Using Theorems \ref{thm:locallimit_excursion} and \ref{thm:ballot1}, we have
\[
 \P\left(S_j \geq -x, S_{\alpha^{(n)}_1} \in [3n^{1/2}, 4n^{1/2}]\right) \geq \frac{c(1 + x)}{n^{1/2}}.
\]
On the other hand, if $1 \not \in F$, for all $h>3$
\begin{multline*}
  \P\left[ S_n \leq y+1, S \in \Upsilon^F(x,y,\delta) \right] \\
  \geq \P\left(S_j \geq -\delta k , \left|S_{\alpha^{(n)}_1}\right| \in [3 n^{1/2}, hn^{1/2}]\right)\qquad \qquad \qquad \qquad \\
  \times \inf_{u \in [3 n^{1/2},h n^{1/2}]}\P_{\alpha^{(n)}_1, u} \left( S_{n-\alpha^{(n)}_1} \leq y+1, S \in \Omega^F(\delta, y) \right).
\end{multline*}
By Theorem \ref{thm:hsurobbins}, we have $\P(\forall n \in \N, S_n \geq - \delta n) > 0$. Thus, writing $\lambda^{(n)} = \floor{\alpha^{(n)}_1/2}$, by central limit theorem, there exists $c>0$ and $h>0$ such that for all $n \geq 1$ large enough
\[
  \P\left(S_j \geq - \delta j, j \leq \lambda^{(n)}, S_{\lambda^{(n)}} \in [-h\sqrt{n}, h \sqrt{n}]\right) \geq c.
\]
Moreover, by Donsker theorem
\[
  \liminf_{n \to +\infty} \inf_{|z| \leq h\sqrt{n}} \P_z\left(S_j \geq -2 h \sqrt{n}, S_{\lambda^{(n)}} \in [3 \sqrt{n}, 4 \sqrt{n}] \right) > 0.
\]
As a consequence, we have
\[P\left(S_j \geq -x \ind{1 \in F} - \delta k \ind{1 \not \in F}, S_{\alpha^{(n)}_1} \in [3n^{1/2}, 4n^{1/2}]\right) \geq \frac{c(1 + x\indset{F}(1))}{n^{\indset{F}(1)/2}}.\]

We now apply time-reversal, for $k \leq n$, let $\hat{S}_k = S_n - S_{n-k}$, we observe that
\begin{multline*}
  \inf_{z \in [3n^{1/2},4n^{1/2}]}\P_{\alpha^{(n)}_1, z} \left( S_{n-\alpha^{(n)}_1} \leq y+1, S \in \Omega^F(\delta, y) \right)\\
  \geq \inf_{u \in [2n^{1/2},5n^{1/2}]} \P\left[\begin{array}{c}\hat{S}_{n-\alpha^{(n)}_1} \in [u,u+1] , \hat{S}_j \geq -\delta n \ind{3 \not \in F}, j \leq n- \alpha^{(n)}_{P-1}\\
  \hat{S}_j \geq n^{1/2}, j \leq n - \alpha^{(n)}_1\end{array}\right].
\end{multline*}
We write $\bar{S}_k = \hat{S}_{n-\alpha^{(n)}_{P-1} + k} - \hat{S}_{n-\alpha^{(n)}_{P-1}}$, we apply again the Markov property at time $n - \alpha^{(n)}_{P-1}$
\begin{multline*}
  \inf_{z \in [3n^{1/2},4n^{1/2}]}\P_{\alpha^{(n)}_1, z} \left( S_{n-\alpha^{(n)}_1} \leq y+1, S \in \Omega^F(\delta, y) \right)\\
  \geq \frac{c}{n^{\indset{F}(3)/2}} \inf_{z \in [0,10n^{1/2}]} \P\left[ \min_{j \leq \alpha^{(n)}_{P-1} - \alpha^{(n)}_1} \bar{S}_j \geq -n^{1/2}, \bar{S}_{\alpha^{(n)}_{P-1} - \alpha^{(n)}_1} \in [z,z+1] \right],
\end{multline*}
using the same tools as above. Finally
\[\P\left[ \min_{j \leq \alpha^{(n)}_{P-1} - \alpha^{(n)}_1} \bar{S}_j \geq -n^{1/2}, \bar{S}_{\alpha^{(n)}_{P-1} - \alpha^{(n)}_1} \in [z,z+1] \right] \geq \frac{c}{n^ {1/2}},\]
using Theorem \ref{thm:locallimit_excursion} and the fact that $\inf_{n \in \N} \P\left[ \min_{j \leq \alpha^{(n)}_{P-1} - \alpha^{(n)}_1} \bar{S}_j \geq -n^{1/2}\right] > 0$, by Donsker's theorem.
\end{proof}

\section{Lagrange multipliers for the optimization problem}
\label{app:optimization_problem}
In this section, for any $\mathbf{h},\mathbf{k} \in \R^P$, we write $\mathbf{h}.\mathbf{k} = \sum_{p=1}^P h_p k_p$ the usual scalar product in $\R^P$. Moreover, if $f : \R^P \to \R$ is differentiable at point $\mathbf{h}$, we write $\nabla f(\mathbf{h}) = \left( \partial_1 f(\mathbf{h}), \ldots \partial_P f(\mathbf{h}) \right)$ the gradient of $f$.

We study in this section the optimization problem consisting of finding $\bba \in \calR$ such that
\begin{equation}
  \label{eqn:optimization}
  \sum_{p=1}^P \left(\alpha_p - \alpha_{p-1}\right)a_p = \sup\left\{ \sum_{p=1}^P \left(\alpha_p-\alpha_{p-1}\right)b_p : \bbb \in \calR \right\}.
\end{equation}
Equation \eqref{eqn:optimization} is a problem of optimization under constraint the $\bba \in \calR$. To obtain a solution, we use an existence of Lagrange multipliers theorem. The version we use here is stated in \cite{Kur76}, for Banach spaces.

\begin{theorem}[Existence of Lagrange multipliers]
\label{thm:lagrange}
Let $P,Q \in \N$. We denote by $U$ an open subset of $\R^P$, $J$ a differentiable function $U \to \R$ and $g=(g_1,\ldots g_Q)$ a differentiable function $U \to \R^Q$. Let $R$ be a convex cone in $\R^Q$ i.e. a subset such that $\forall x, y \in R, \forall \lambda, \mu \in \R_+, \lambda x + \mu y \in R$.

If $\bba \in \R^P$ verifies $g(a) \in R$ and
\[
  J(\bba) = \sup\left\{ J(\bbb), \bbb \in \R^p : g(\bbb) \in R \right\},
\]
and if the differential of $g$ at point $\bba$ is a surjection, then there exist non-negative \textit{Lagrange multipliers} $\lambda_1,\ldots \lambda_Q$ verifying the following properties.
\begin{description}
  \item[(L1)] For all $\mathbf{h} \in \R^P$, $\nabla J(\bba) . \mathbf{h} = \sum_{q=1}^Q \lambda_q (\nabla g_q(a) . \mathbf{h})$.
  \item[(L2)] For all $h \in R$, $ \sum_{q=1}^Q \lambda_q h_q \leq 0$;
  \item[(L3)] $\sum_{q=1}^Q \lambda_q g_q(\bba) = 0$. 
\end{description}
\end{theorem}

Using this theorem, we prove Proposition \ref{prop:optimization_problem}. We start by proving that if $\bba$ satisfies some specific properties, then $\bba$ is the solution to \eqref{eqn:optimization}.

\begin{lemma}
\label{lem:conditions}
Under assumptions \eqref{eqn:breeding} and \eqref{eqn:differentiable}, $\bba \in \calR$ is a solution of \eqref{eqn:optimization} if and only if, writing $\theta_p = \left( \kappa^*_p \right)'(a_p)$, we have
\begin{description}
  \item[(P1)] $ \bbtheta$ is non-decreasing and positive ;
  \item[(P2)] if $K^*(\bba)_p<0$, then $\theta_{p+1}=\theta_p$ ;
  \item[(P3)] $K^*(\bba)_P=0$.
\end{description}
\end{lemma}

\begin{proof}
For $\bbb \in \R^P$, we denote by
\[ J(\bbb) = \sum_{p=1}^P (\alpha_p-\alpha_{p-1})b_p, \quad R = \{ \mathbf{k} \in \R^P : k_p \leq 0, p \leq P\}, \]
and we write, $\theta_p(\bbb) = (\kappa^*_p)'(b_p)$.

We assume in a first time that $\bba \in \calR$ is a solution of \eqref{eqn:optimization}, in which case
\begin{equation}
  \label{eqn:opt}
  J(\bba) = \sup\left\{ J(\bbb), \bbb \in \R^P : K^*(\bbb) \in R \right\}.
\end{equation}
The function $J$ is linear thus differentiable, and assumption \eqref{eqn:differentiable} implies that $K^*$ is differentiable at point $\bba$. For $\mathbf{h} \in \R^P$, we have $\nabla J(\bba) . \mathbf{h} = \sum_{p=1}^P (\alpha_p - \alpha_{p-1})h_p,$ and $\nabla K^*(\bba)_p . \mathbf{h}  =  (\alpha_p-\alpha_{p-1})\theta_p(\bba) h_p$.

To prove that $K^*$ has a surjective differential, it is enough to prove that for all $p \leq P$, $\theta_p(\bba) \neq 0$. Let $p \leq P$ be the smallest value such that $\theta_p(\bba)=0$. Observe that in this case, $\kappa^*_p(a_p)<0$ by \eqref{eqn:breeding}, thus we can increase a little $a_p$ and stay in $\calR$ as soon as we decrease a little $a_{p-1}$ --or $a_P$ if $p=1$, in which case same proof would work with few modifications. For $\epsilon>0$ and $q \leq P$, we write $\bba^\epsilon_q = \bba_q - \epsilon \ind{q=p-1}+ \epsilon^{2/3} \ind{q=p}$. We observe that, for all $\epsilon>0$ small enough,
\begin{align*}
  K^*(\bba^\epsilon)_{p-1} & = K^*(\bba^\epsilon)_{p-2} + (\alpha_{p-1} - \alpha_{p-2}) \kappa^*_{p-1} (a_{p-1} - \epsilon)\\
  &\leq K^*(\bba)_{p-2} + (\alpha_{p-1}-\alpha_{p-2}) \kappa^*_{p-1}(a_{p-1}) - (\alpha_{p-1} - \alpha_{p-2})\theta_{p-1}(\bba) \epsilon + O(\epsilon^2)\\
  &\leq K^*(\bba)_{p-1} - (\alpha_{p-1} - \alpha_{p-2})\theta_{p-1}(\bba) \epsilon + O(\epsilon^2)
\end{align*}
and
\begin{align*}  
  K^*(\bba^\epsilon)_p &\leq  K^*(\bba)_{p-1} + (\alpha_{p}-\alpha_{p-1})\kappa^*_p(a_p + \epsilon^{2/3})\\
  &\leq  K^*(\bba)_{p-1} - (\alpha_{p-1} - \alpha_{p-2})\theta_{p-1}(\bba) \epsilon + (\alpha_p-\alpha_{p-1}) \kappa^*_p(a_p)+ O(\epsilon^{4/3})\\
  &\leq K^*(\bba)_p - (\alpha_{p-1} - \alpha_{p-2})\theta_{p-1}(\bba) \epsilon+ O(\epsilon^{4/3}),
\end{align*}
thus, for $\epsilon >0$ small enough, $\bba^\epsilon \in \calR$ and $\sum_{p=1}^P (\alpha_p -\alpha_{p-1})a^\epsilon_p > \sum_{p=1}^P (\alpha_p -\alpha_{p-1})a_p$, which is inconsistent with the fact that $\bba$ is the optimal solution of \eqref{eqn:optimization}.

Therefore, by Theorem \ref{thm:lagrange}, there exist non-negative $\lambda_1,\ldots \lambda_P$ such that
\begin{description}
  \item[(L1)] $\forall \mathbf{h} \in \R^P$, $\nabla J(a). \mathbf{h} = \sum_{p=1}^P \lambda_p \nabla K^*(\bba)_p . \mathbf{h}$;
  \item[(L2)] $\forall \mathbf{h} \in R$, $\sum_{p=1}^P \lambda_p h_p \leq 0$;
  \item[(L3)] $\sum_{p=1}^P \lambda_p K^*(\bba)_p = 0$. 
\end{description}

We observe that Condition (L1) can be rewritten $\forall p \leq P, \lambda_p \theta_p(\bba) = 1$, therefore $\theta_p(\bba) = \frac{1}{\lambda_p}$. Moreover, Condition (L2) applied to the vector $\mathbf{h}^p \in R$ defined by $h^p_j = -\ind{j=p} + \ind{j=p+1}$ implies that $\lambda$ is non-increasing, thus $\bbtheta$ is non-decreasing; which gives (P1). Finally, we rewrite Condition (L3) as follows, by discrete integration by part
\[
  0 = \sum_{p=1}^P \lambda_p K^*(\bba)_p = \underbrace{\lambda_P K^*(\bba)_P}_{\leq 0} - \sum_{p=1}^{P-1} \underbrace{(\lambda_{p+1}-\lambda_p)K^*(\bba)_p}_{\geq 0},
\]
therefore Condition (P3) ($K^*(\bba)_P=0$) is verified; and if $\lambda_{p+1}\neq \lambda_p$, then $K^*(\bba)_p=0$, which implies (P2).

We now suppose that $\bba \in \calR$ verifies Conditions (P1), (P2) and (P3) and we prove that for all $\bbb \in \calR$,
\begin{equation}
  \label{eqn:inter}
  \sum_{p=1}^P (\alpha_p- \alpha_{p-1}) a_p \geq \sum_{p=1}^P (\alpha_p - \alpha_{p-1}) b_p.
\end{equation}
To do so, we use the fact that functions $\kappa^*_p$ are convex and differentiable at point $\bba$, therefore, for all $x \in \R$, $\kappa^*_p(x) \geq \kappa^*_p(a_p) + \theta_p (x-a_p)$. As a consequence, we have
\begin{align*}
  \sum_{p=1}^P (\alpha_p- \alpha_{p-1}) (a_p-b_p)
  &\geq \sum_{p=1}^P \frac{\kappa^*_p(a_p) - \kappa^*_p(b_p)}{\theta_p} (\alpha_p - \alpha_{p-1})\\
  &\geq (K^*(\bba)_P - K^*(\bbb)_P) \frac{1}{\theta_P} - \sum_{p=1}^{P-1} \left( \frac{1}{\theta_{p+1}} - \frac{1}{\theta_p} \right) (K^*(\bba)_p - K^*(\bbb)_p)
\end{align*}
by discrete integration by part. By the specific properties of $a$, we have
\[K^*(\bba)_P \frac{1}{\theta_P} - \sum_{p=1}^{P-1}\left( \frac{1}{\theta_{p+1}} - \frac{1}{\theta_p} \right) K^*(\bba)_p = 0,\]
thus
$\displaystyle 
  \sum_{p=1}^P (\alpha_p- \alpha_{p-1}) (a_p-b_p) \geq - \frac{K^*(\bbb)_P}{\theta_P} + \sum_{p=1}^{P-1} \left( \frac{1}{\theta_{p+1}} - \frac{1}{\theta_p} \right) K^*(\bbb)_p \geq 0,
$
as $\bbtheta$ is non-decreasing an $K^*(\bbb)$ non-positive. Optimizing \eqref{eqn:inter} over $\bbb \in \calR$ gives us
\[
  \sum_{p=1}^P (\alpha_p - \alpha_{p-1})a_p \geq v_{\mathrm{is}}
\]
which ends the proof.
\end{proof}

We now prove the uniqueness of the solution of \eqref{eqn:optimization}.
\begin{lemma}
If for all $p \leq P$, $\kappa_p$ is finite on an open subset of $[0,+\infty)$, then there is at most one solution to \eqref{eqn:optimization}.
\end{lemma}

\begin{proof}
The uniqueness of the solution in an easy consequence of the strict convexity of $(\kappa^*_p,p\leq P)$. Let $\bba$ and $\bbb$ be two different solutions to \eqref{eqn:optimization}, there exists a largest $p \leq P$ such that $a_p \neq b_p$. Writing $\mathbf{c}=\frac{\bba+\bbb}{2}$, for any $q \geq p$, we have $K^*(\mathbf{c})_q < \frac{K^*(\bba)_q + K^*(\bbb)_q}{2} \leq 0$. Thus, by continuity of $K^*$, $\mathbf{c}$ is in the interior of $\calR$, then we can increase a little $c_p$, and the path driven by $(\mathbf{c} + \epsilon \ind{.=p})$ goes farther than both $\bba$ and $\bbb$, which is a contradiction.
\end{proof}

Finally, we prove the existence of such a solution when the mean number of children of an individual in the BRWis is finite.
\begin{lemma}
\label{lem:existence}
Under the assumptions \eqref{eqn:differentiable} and \eqref{eqn:finite_reproduction}, there exists at least a solution to \eqref{eqn:optimization}.
\end{lemma}

\begin{proof}
If $\kappa_p(0)<+\infty$, then $\inf_{\R} \kappa^*_p = - \kappa_p(0)$ and the minimum is reached at $\kappa'_p(0)$. As $\kappa^*_p$ are bounded from below, for all $p \leq P$ there exists $x_p \geq 0$ such that
\[
   (\alpha_p - \alpha_{p-1}) \kappa^*_p(x_p) + \sum_{q \neq p} (\alpha_q - \alpha_{q-1}) \inf_{\R} \kappa^*_q > 0.
\]
Therefore, writing $X = \calR \cap \prod_{p \leq P} [\kappa'_p(0), x_p ]$, we have
\[
  \sup_{\bbb \in \calR} \sum (\alpha_p - \alpha_{p-1}) b_p = \sup_{\bbb \in X} \sum (\alpha_p - \alpha_{p-1}) b_p.
\]
But, $X$ being compact, this supremum is in fact a maximum. There exists $\bba \in X$ such that $\displaystyle \sum (\alpha_p - \alpha_{p-1}) a_p =  \sup_{\bbb \in \calR} \sum (\alpha_p - \alpha_{p-1}) b_p$
which ends the proof.
\end{proof}

\section{Notation}

\begin{itemize}
  \item {\em Point processes}
  \begin{itemize}
    \item $\calL_p$: law of a point process;
    \item $L_p$: point process with law $\calL_p$;
    \item $\kappa_p$: log-Laplace transform of $\calL_p$;
    \item $\kappa^*_p$: Fenchel-Legendre transform of $\calL_p$;
    \item $X_p$: defined in \ref{eqn:definexp};
    \item $v_p = \inf_{\theta > 0} \frac{\kappa_p(\theta)}{\theta}$: speed of branching random walk with reproduction law $\calL_p$;
    \item $\bar{\theta}_p$ critical parameter such that $\bar{\theta}_p v_p - \kappa_p(\bar{\theta}_p)=0$;
  \end{itemize}
  \item {\em Generic marked tree}
  \begin{itemize}
    \item $\T$: genealogical tree of the process;
    \item $u \in \T$: individual in the process;
    \item $V(u)$: position of the individual $u$;
    \item $|u|$: generation at which $u$ belongs;
    \item $u_k$: ancestor at generation $k$ of $u$;
    \item $\emptyset$: initial ancestor of the process;
    \item if $u \neq \emptyset$, $\pi u$: parent of $u$;
    \item $\Omega(u)$: set of the children of $u$;
    \item $M_n = \max_{|u|=n} V(u)$ maximal displacement at the $n^\text{th}$ generation in $(\T,V)$.
  \end{itemize}
  \item {\em Branching random walk through a series of interfaces}
  \begin{itemize}
    \item $P$: number of distinct phases in the process;
    \item $0=\alpha_0 < \alpha_1 < \ldots < \alpha_P=1$: position of the interfaces;
    \item $\alpha^{(n)}_p = \floor{n \alpha_p}$: position of the $p^\text{th}$ interface for the BRWis of length $n$;
    \item $\bar{a}^{(n)}_k=\sum_{p=1}^P \sum_{j=1}^k \ind{j \in (\alpha^{(n)}_{p-1},\alpha^{(n)}_p]}$ path driven by $\bba := (a_1,\ldots a_p) \in \R^P$;
    \item $u$ ``follows path $\bar{a}^{(n)}$'' if $\forall k \leq |u|$, $|V(u_k)-\bar{a}^{(n)}_k| \leq n^{1/2}$;
    \item $K^*(\bba)_p =  \sum_{q=1}^p (\alpha_q - \alpha_{q-1})a_q$: rate function associated to the BRWis;
    \item $\calR = \left\{ \bba \in \R^P : \forall p \leq P, K^*(\bba)_p \leq 0 \right\}$: set of $\bba \in \R^P$ such that $\bar{a}^{(n)}$ is followed until time $n$ by at least one individual with positive probability.
  \end{itemize}
  \item {\em The optimal path}
  \begin{itemize}
    \item $v_\mathrm{is} = \max_{\bbb \in \calR} \sum_{p=1}^P (\alpha_p-\alpha_{p-1})b_p$: speed of the BRWis;
    \item $\bba \in \calR$ such that $\sum_{p=1}^P (\alpha_p-\alpha_{p-1}) a_p = v_\mathrm{is}$: optimal speed profile;
    \item $\theta_p = (\kappa^*_p)'(a_p)$;
    \item $T = \# \{ \theta_p, p \leq P\}$: number of different values taken by $\bbtheta$;
    \item $\phi_1<\phi_2 < \cdots < \phi_T$: different values taken by $\bbtheta$;
    \item $f_t = \min\{k \leq P : \theta_k = \phi_t\}$ and $l_t=\max\{k \leq P : \theta_k = \phi_t\}$;
    \item $\lambda= \sum_{t=1}^T \frac{1}{2\phi_t}\left[ \ind{K^*(\bba)_{f_t}=0} + 1 + \ind{K^*(\bba)_{l_t-1}=0} \right]$: logarithmic correction;
    \item $B = \{p \leq P : K^*(\bba)_{p-1} = K^*(\bba)_p = 0 \}$: phases such that the optimal path is close to the boundary of the BRWis;
  \end{itemize}
  \item {\em Spinal decomposition}
  \begin{itemize}
    \item $W_n = \sum_{|u|=n} e^{\theta V(u) - \sum_{k=1}^n \kappa_k(\theta)}$: the additive martingale with parameter $\theta$;
    \item $\P_{k,x}$: law of the time-inhomogeneous branching random walk with environment $(\calL_k,\calL_{k+1},\ldots)$;
    \item $\bar{\P}_{k,x} = W_n \cdot \P_{k,x}$: size-biased law of $\P_{k,x}$;
    \item $\hat{\P}_{k,x}$: law of the branching random walk with spine;
    \item $w$: spine of the branching random walk;
    \item $\calF_n = \sigma(u,V(u), |u| \leq n)$: filtration of the branching random walk;
    \item $\calG_n = \sigma(w_k, V(w_k), k \leq n) \vee \sigma(u,V(u), u \in \Omega(w_k), k < n)$: filtration of the spine;
    \item $\hat{\calF}_n = \calF_n \vee \calG_n$: filtration of the branching random walk with spine;
    \item Spinal decomposition: Proposition \ref{prop:spinaldecomposition};
    \item Many-to-one lemma: Lemma \ref{lem:manytoone}.
  \end{itemize}
  \item {\em Random walks}
  \begin{itemize}
    \item $(T_n)$: random walk with finite variance;
    \item $(S_n)$: random walk through a series of interfaces, its law under $\P_{k,x}$ is the same as the law of $(V(w_j), j \leq n-k)$ under $\hat{\P}_{k,x}$.
    \item Time-reversal: replace $S$ by the random walk $(\hat{S}_n = S_n - S_{n-k}, k \leq n)$.
  \end{itemize}
  \item {\em Branching random walk estimates}
  \begin{itemize}
    \item $m_n = nv_\mathrm{is} - \lambda \log n$;
    \item $E_p(\phi) = \sum_{q=1}^p (\alpha_q-\alpha_{q-1})(\phi \kappa'_p(\phi) - \kappa_q(\phi))$;
    \item $K^{(n)}_k = \sum_{p=1}^P \kappa_p(\theta) \sum_{j=1}^k \ind{j \in (\alpha^{(n)}_{p-1},\alpha^{(n)}_p]}$;
    \item $r^{(n)}_k = a_P(k-n) + \frac{3}{2\theta} \log (n-k+1)$;
    \item $B^{(n)} = \bigcup_{p \in B} (\alpha^{(n)}_{p-1}, \alpha^{(n)}_p]$ and $F^{(n)} = \bigcup_{p \in B \cap \{1,P\}} [\alpha^{(n)}_{p-1}, \alpha^{(n)}_p]$;
    \item $f^{(n)}_j = a_1 j \ind{j \leq \alpha^{(n)}_1} + \left(m_n + r^{(n)}_k\right) \ind{j \geq \alpha^{(n)}_{P-1}}$;
    \item $X^{(n)}(y,h) = \sum_{|u|=n} \ind{V(u) - m_n -y \in [-h,-h+1]} \ind{V(u_j) \leq f^{(n)}_j + y, j \in F^{(n)}}$;
    \item for $\delta>0$ such that $3\theta \delta < \min_{p \in B^c} -E_p(\theta)$,
    \[ g^{(n)}_k = 1+
\begin{cases}
  \bar{a}^{(n)}_k - \ind{p=P} \lambda \log n & \mathrm{if} \quad E_p(\theta)=E_{p-1}(\theta)=0\\
  \bar{a}^{(n)}_k + (k-\alpha^{(n)}_{p-1}) \delta & \mathrm{if} \quad E_{p-1}(\theta)=0,E_p(\theta)<0\\
  \bar{a}^{(n)}_k + (\alpha^{(n)}_p-k) \delta & \mathrm{if} \quad E_p(\theta)=0, E_{p-1}(\theta)<0\\
  \bar{a}^{(n)}_k + \delta n &\mathrm{otherwise;}
\end{cases} \]
    \item $\calA_n(y) = \left\{ |u| = n : V(u) \geq m_n + y, V(u_j) \leq g^{(n)}_j + y, j \leq n \right\}$;
    \item $\xi(u) = \sum_{u' \in \Omega(u)} (1 + (V(u')-V(u))_+\ind{|u| \in B^{(n)}})e^{\theta (V(u')-V(u))}$;
    \item $\calB_n(z) = \left\{ |u|=n : \xi(u_j) \leq z e^{-\tfrac{\theta}{2} \left[ V(u_j) - g^{(n)}_j \right]} \right\}$;
    \item $G_n(y,z) = \calA_n(y) \cap \calB_n(z)$ and $Y_n(y,z) = \# G_n(y,z)$
  \end{itemize}
\end{itemize}

\textbf{Acknowledgements.} I wish to thank my supervisor Zhan Shi for his constant help while working on this subject, all the useful discussions and advices, Ofer Zeitouni for his explanations on \cite{FaZ12a}, and the anonymous referee for his helpful comments on the earlier versions of this article. I also wish to thank Ming Fang for pointing out a mistake in an earlier version of this article, and Denis Denisov for pointing \cite{DSV} to me.

\nocite{*}

\bibliographystyle{plain}

\end{document}